\documentclass[12pt,a4paper]{article}

\pdfoutput=1
\usepackage{amsmath}
\usepackage{amsfonts}
\usepackage{amssymb}
\usepackage{amsthm}
\usepackage{graphicx}
\usepackage[all]{xy}

\newtheorem{thm}{Theorem}[section]
\newtheorem{prop}{Proposition}[section]

\newtheorem{dfn}{Definition}[section]

\theoremstyle{remark}
\newtheorem{rmk}{Remark}[section]

\newtheorem{notation}{Notation}[section]

\newcommand{\R}{\mathbb{R}}
\newcommand{\funcman}[1]{C^{\infty}(#1)}
\newcommand{\manifold}[1]{\mathbf{#1}}
\newcommand{\atlas}[2]{\mathfrak{#1}(\manifold{#2})}
\newcommand{\categories}[1]{\mathbf{#1}}
\newcommand{\liebr}[2]{[{#1},{#2}]}

\newcommand{\field}[1]{\mathbb{#1}}
\newcommand{\catvect}[1]{\categories{Vec}^{\field{#1}}}

\author{Luiz Henrique Pereira P\^egas\footnote{Contact by e-mail: lhp@ime.usp.br or divtzero@gmail.com}}
\title{About the Hochschild-Kostant-Rosenberg theorem for differentiable manifolds}
\date{\today}

\begin{document}
\maketitle

\begin{abstract}
In this notes it will be provided a set of techniques which can help one to understand the proof of the Hochschild-Kostant-Rosenberg theorem for differentiable manifolds. Precise definitions of multidiferential operators and polyderivations on an algebra are given, allowing to work on these concepts, when the algebra is an algebra of functions on a differentiable manifold, in a coordinate free description. Also, it will be constructed a cup product on polyderivations which corresponds on (Hochschild) cohomology to wedge product on multivector fields. At the end, a proof of the above mentioned theorem will be given.
\end{abstract}

\tableofcontents

\section*{Acknowledgements}
I wish to thanks prof. Dr. Eduardo Hoefel, my advisor on my MSc dissertation in which this notes are based, for the uncountable hours of discussion and guidance. Without his patience, this work could not be possible.

\section{Introduction}
The main purpose of this notes is to provide a set of techniques which can help one to understand the proof of the Hochschild-Kostant-Rosenberg for differentiable manifolds. On the way to do that, precise definitions of multidiferential operators and polyderivations on an algebra are given. When the algebra is an algebra of functions on a differentiable manifold, this allows us to work on these concepts in a coordinate free description. Also, it is constructed a cup product on polyderivations which corresponds on (Hochschild) cohomology to wedge product on multivector fields.

To fully understand this notes, some background on algebra and differentiable manifolds is desirable. In section 2 some basic concepts are presented and some notation are fixed. So, the readers who have some knowledge on differential geometry, derivations and differential operators on associative algebras may goes directly to section 3. In section 3 the concept of multiderivations on an associative, commutative, unital algebra is defined and some results on algebras of smooth functions on a manifold are stated. In section 4 the concept of iterated derivations is defined and related to derivations of higher orders on the algebra of smooth functions on a manifold. In section 5 the concept of polyderivations of an associative, commutative, unital algebra is defined and it is used to provide a coordinate free version of polyderivantions on a manifold. Section 5 ends with a proof of the Hochschild-Kostant-Rosenberg theorem for differentiable manifolds which uses the tools introduced in the previous sections.

Throughout this notes, $Hom_{\mathbf{C}}(A,B)$ will denote the morphisms from $A$ to $B$ in the category $\mathbf{C}$ and $A \approx_{\mathbf{C}} B$ will mean that $A$ is isomorphic to $B$ in the category $\mathbf{C}$. As an example, $Hom_{\catvect{K}}(A,B)$ is the set of $\field{K}$-linear transformations between the $\field{K}$-vector spaces $A$ and $B$.
\section{Algebras, Derivations and the Tangent Bundle}
\label{algebras}
In this section it will be provided some of the standard concepts needed to understand the techniques developed in the final sections of this notes. It was made for setting some notations and some ``ways of thinking''. From now on, $k$ will denote a commutative unital ring and $\field{K}$ will denote a field.
\begin{dfn}[Derivation on an algebra]
Let $(A,\mu,e)$ be an associative $k$-algebra with unity $e$ (briefly $A$, when the product $\mu$ is clear from the context). A derivation $D$ on $A$ is a $k$-linear map $D:A \rightarrow A$ such that
\begin{displaymath}
D(\mu)=\mu(D \otimes id)+\mu(id \otimes D)
\end{displaymath}
where $id$ is the identity on $A$.\label{derivation_algebra}
\end{dfn}
The condition above is known as ``Leibniz rule''.
\begin{rmk}
The unity $e$ on $A$ provides an immersion $e:k \rightarrow A$, so elements in $k$ can be viewed as elements in $A$ through such immersion. If $D$ is a derivation on $A$, by the Leibniz rule we have for all $a \in k$
\begin{displaymath}
D(a)=aD(e)=aD(\mu(e \otimes e))=a(\mu(D \otimes e)+\mu(e \otimes D))=D(a)+D(a)
\end{displaymath}
therefore $D(a)=0$.
\end{rmk}
\begin{notation}
The $k$-module of all derivations on $A$ will be denoted by $Der(A)$.
\end{notation}
\begin{dfn}[Inner derivations]
Let $(A,\mu,e)$ be an associative $k$-algebra with unity $e$. A $k$-linear map $f:A \rightarrow A$ is an inner derivation on $A$ if and only if there exists $a \in A$ such that
\begin{displaymath}
f=\mu(a \otimes id)-\mu(id \otimes a),
\end{displaymath}
where $a \otimes id$ denotes $a \otimes id:A \rightarrow A \otimes A$ such that $(a \otimes id)(b) = a \otimes b$.\label{inner_derivation}
\end{dfn}
\begin{rmk}
By denoting $IDer(A)$ the $k$-module of inner derivations on $A$, we have $IDer(A) \subset Der(A)$.
\end{rmk}
Now we make precise the notion of high order differential operator and high order derivation on a commutative associative $k$-algebra. The definition, naturally, is recursive.
\begin{dfn}[Higher order differential operator]
Let $(A,\mu,e)$ be a commutative associative $k$-algebra with unity $e$. A $k$-linear map $D: A \rightarrow A$ is a differential operator of order $\leq$ $r$ on $A$, with $r \in \mathbb{N} \setminus \{0\}$ if and only if for all $a \in A$ the map
\begin{displaymath}
\tilde{D}=D(\mu(a \otimes id))-\mu(a \otimes D)
\end{displaymath}
is a differential operator of order $\leq$ $r-1$. A map $D: A \rightarrow A$ is a differential operator of order 0 if and only if there exists $a \in A$ such that $D=\mu(a \otimes id)$.\label{high_order_diff_operator}
\end{dfn}
\begin{dfn}[Higher order derivation]
Let $(A,\mu,e)$ be a commutative associative $k$-algebra with unity $e$. A derivation of order $\leq$ $r$ on $A$ is a differential operator of order $\leq$ $r$ $D$ such that $D(a)=0,\;\forall a \in k$.\label{high_order_derivation}
\end{dfn}
\begin{thm}
Let $(A,\mu,e)$ be a commutative associative $k$-algebra with unity $e$. Then derivations on $A$ are exactly the derivations of order $\leq$ 1 on $A$ and a differential operator of order $\leq$ 1 $D$ can be written uniquely as $D = \partial + D(e)$ with $\partial \in Der(A)$.
\begin{proof}
Here it is convenient to use de juxtaposition to denote the product $\mu$. If $D \in Der(A)$ then $D$ is a derivation of order $\leq$ 1 because given $a \in A$
\begin{displaymath}
\tilde{D}(b)=D(ab)-aD(b)=D(a)b
\end{displaymath}
and if $a \in k$ we have $D(a)=0$.

Conversely, if $D:A \rightarrow A$ is a derivation of order $\leq$ 1 then the map $\partial = D - D(e)$, \emph{i.e.}, $\partial(a) = D(a) - aD(e),\;\forall a \in A$ satisfies
\begin{displaymath}
\partial(a)=D(a)-aD(e)=f_{a}
\end{displaymath}
were $f_{a}$ is the element of $A$ given by
\begin{displaymath}
D(ab) - aD(b) = f_{a}b
\end{displaymath}
(just do $b=e$) thus
\begin{displaymath}
D(ab)-aD(b)=\partial(a)b
\end{displaymath}
hence
\begin{eqnarray*}
\partial(ab) & = & D(ab) - abD(e) = D(ab) - aD(b) + aD(b) - abD(e) = \\
& = & \partial(a)b + aD(b)-abD(e) = \partial(a)b + a(D(b) - bD(e)) = \\
& = & \partial(a)b + a\partial(b)
\end{eqnarray*}
which means $\partial \in Der(A)$.
\end{proof}
\end{thm}
When the algebra $(A,\mu)$ is graded we can define a derivation which ``respects'' such structure.
\begin{dfn}[Graded derivation]
Let $(A,\mu)$ be a graded $k$-algebra. A $k$-linear map $D:A \rightarrow A$ is a graded derivation on $A$ of degree $p$ if and only if for all homogeneous elements $a \in A_{i}$ and for all $b \in A$ it satisfies
\begin{displaymath}
D(\mu(a \otimes b)) = \mu(D(a) \otimes b)+(-1)^{pi}\mu(a \otimes D(b))
\end{displaymath}
\end{dfn}
Now it is convenient to set some notions (and notations).
\begin{dfn}[Superalgebra]
We say that a graded $k$-algebra $(A,\mu)$ is a superalgebra (or supercommutative) if and only if $\mu$ satisfies for all $a \in A_{i}$, $b \in A_{j}$
\begin{displaymath}
\mu(a \otimes b)=(-1)^{ij}\mu(b \otimes a)
\end{displaymath}
\end{dfn}
\begin{dfn}[Lie superalgebra]
A Lie superalgebra is a pair $(L,\liebr{\;}{\;})$ where $L$ is a $\mathbb{Z}$-graded $\field{K}$-vector space and $\liebr{\;}{\;}:L \times L \rightarrow L$ is a bilinear map such that
\begin{enumerate}
\item[i)] $\liebr{L_{i}}{L_{j}} \subset L_{i+j},\; \forall i,j \in \mathbb{Z}$;
\item[ii)] $\liebr{a}{b}=-(-1)^{ij}\liebr{b}{a},\;\forall a \in L_{i},\; \forall b \in L_{j}$,  (graded antisymmetry);
\item[iii)] $\liebr{a}{\liebr{b}{c}}=\liebr{\liebr{a}{b}}{c}+(-1)^{ij}\liebr{b}{\liebr{a}{c}},\;\forall a \in L_{i},\; \forall b \in L_{j},\; \forall c \in L$, (graded Jacobi).
\end{enumerate}
\end{dfn}
\begin{dfn}[Derivation on a Lie superalgebra]
Let $(L,\liebr{\;}{\;})$ be a Lie superalgebra. A $\field{K}$-linear map $D:L \rightarrow L$ is a derivation of degree $p$ on $L$ if and only if it satisfies for all $a \in L_{i}$ and for all $b \in L$
\begin{displaymath}
D(\liebr{a}{b})=\liebr{D(a)}{b}+(-1)^{ip}\liebr{a}{D(b)}
\end{displaymath}
\end{dfn}
\begin{rmk}
It is easy to see that if $(L,\liebr{\;}{\;})$ is a Lie superalgebra, then for any $a \in L_{i}$, $ad(a):L \rightarrow L$ given by $ad(a)(b)=\liebr{a}{b}$ is a degree $i$ derivation on $L$. This is exactly what the Jacobi identity is about.
\end{rmk}

In differential geometry, the concept of tangent vector on a manifold at a point is often given by using equivalence classes of curves or simply stating the property of being a derivation at a point. In this notes we will use an equivalent (for finite dimensional differentiable manifolds) and well known definition for tangent vectors which is slightly different from the usual definition, but reveals some interesting aspects. The construction given here follows \cite{warner_dif_man}.

Let $\manifold{M}$ be a differentiable manifold. For each $p \in \manifold{M}$, define the $\R$-vector space $V_{p}$ tangent at $p$ as the following. Define the relation $\sim$ on $\funcman{\manifold{M}}$ by $f \sim g$ if and only if there exists an open neighbourhood $U$ of $p$ such that $f|_{U}=g|_{U}$. This is an equivalence relation (the reader is invited to prove this, it is not hard) and the equivalence classes induced are called germs of functions at $p$. The algebraic operations on $\funcman{\manifold{M}}$ can be used to induce a $\R$-algebra structure on $\funcman{\manifold{M}} / \sim$. We denote such structure by $\mathcal{F}_{p}$. Let $I_{p}$ be the ideal of $\mathcal{F}_{p}$ of germs of functions that vanish at $p$. As $I_{p}$ is an ideal of $\mathcal{F}_{p}$ and $I_{p}^{2}$ is an ideal of $I_{p}$, we have that $I_{p}/I_{p}^{2}$ is a $\R$-vector space. Then we define $V_{p}=(I_{p}/I_{p}^{2})^{\ast}$, \emph{i.e.}, $V_{p}$ is the vector space dual to $I_{p}/I_{p}^{2}$. We will prove that $V_{p}$ is finite dimensional. From now on, we will denote by $(U,\varphi,m)$ a local chart on a differentiable manifold $\manifold{M}$ where $U$ is an open subset of $\manifold{M}$, $\varphi:U \rightarrow U_{0}$ is a diffeomorphism from $U$ to an open subset $U_{0}$ of $\R^{m}$ or simply $(U,\varphi)$ if the dimension of $\manifold{M}$ is clear.
\begin{prop}
Let $(U,\varphi,m)$ be a local chart of $\manifold{M}$ around $p$. Denoting the $i$-th canonical projection on $\R^{m}$ by $t^{i}:\R^{m} \rightarrow \R$ and the $i$-th coordinate function on $U$ by $x^{i}=t^{i} \circ \varphi$ then the set of equivalence classes of $x^{i}$ for $i=1,\ldots,m$ in $I_{p}/I_{p}^{2}$ constitutes a basis for such space.
\begin{proof}
Given $\mathbf{f} \in I_{p}/I_{p}^{2}$, let $f \in \funcman{\manifold{M}}$ be a representing element. Note that $f(p)=0$. Without loss of generality we can suppose $\varphi(U)$ convex\footnote{It is necessary because we want to use Taylor's formula with integral reminder.} and $\varphi(p)=0$. The coordinate expression of $f$ is given by $f \circ \varphi^{-1}:\varphi(U) \rightarrow \R$, which by the Taylor formula gives, for $y=\varphi(q), q \in U$
\begin{eqnarray*}
(f \circ \varphi^{-1})(y) & = & \sum_{i=1}^{m}\frac{\partial(f \circ \varphi^{-1})}{\partial t^{i}} \biggl |_{0} t^{i}(y) + \\
& + & \sum_{i,j=1}^{m} t^{i}(y)t^{j}(y) \int_{0}^{1}(1-s)\frac{\partial^{2}(f \circ \varphi^{-1})}{\partial t^{i} \partial t^{j}} \biggl |_{sy}ds \quad \\
(f \circ \varphi^{-1})(\varphi(q)) & = & \sum_{i=1}^{m}\frac{\partial(f \circ \varphi^{-1})}{\partial t^{i}} \biggl |_{0} t^{i}(\varphi(q)) + \\
& + & \sum_{i,j=1}^{m} t^{i}(\varphi(q))t^{j}(\varphi(q)) \int_{0}^{1}(1-s)\frac{\partial^{2}(f \circ \varphi^{-1})}{\partial t^{i} \partial t^{j}} \biggl |_{sy}ds \quad \\
f(q) & = & \sum_{i=1}^{m}\frac{\partial(f \circ \varphi^{-1})}{\partial t^{i}} \biggl |_{\varphi(p)} x^{i}(q) + \\
& + & \sum_{i,j=1}^{m} x^{i}(q)x^{j}(q) \int_{0}^{1}(1-s)\frac{\partial^{2}(f \circ \varphi^{-1})}{\partial t^{i} \partial t^{j}} \biggl |_{sy}ds
\end{eqnarray*}
Because $f \in \funcman{\manifold{M}}$ and $x^{i}(p)=0$, the term
\begin{equation*}
\sum_{i,j=1}^{m} x^{i}x^{j} \int_{0}^{1}(1-s)\frac{\partial^{2}(f \circ \varphi^{-1})}{\partial t^{i} \partial t^{j}} \biggl |_{sy}ds
\end{equation*}
represents the null class in $I_{p}/I_{p}^{2}$. From this we infer that we can write
\begin{equation*}
\mathbf{f}=\sum_{i=1}^{m}\frac{\partial(f \circ \varphi^{-1})}{\partial t^{i}} \biggl |_{\varphi(p)} \mathbf{x}^{i}
\end{equation*}
where $\mathbf{x}^{i}$ is the class of $x^{i}$ in $I_{p}/I_{p}^{2}$. Hence the set $\{\mathbf{x}^{i}\}, i=1, \ldots m$ spans $I_{p}/I_{p}^{2}$. To show the linear independence, notice that
\begin{equation*}
\sum_{i=1}^{m}a_{i}\mathbf{x}^{i}=0 \Rightarrow \sum_{i=1}^{m}a_{i}[x^{i}] \in I_{p}^{2}
\end{equation*}
where $[x^{i}]$ is a representing of $\mathbf{x}^{i}$ in $I_{p}$. Writing in coordinates, we have:
\begin{equation*}
\left ( \sum_{i=1}^{m}a_{i}x^{i} \right ) \circ \varphi^{-1} = \sum_{i=1}^{m}a_{i}(x^{i} \circ \varphi^{-1}) = \sum_{i=1}^{m}a_{i}t^{i}
\end{equation*}
which shows that $\displaystyle \sum_{i=1}^{m}a_{i}[t^{i}] \in I_{\varphi(p)}^{2}$, because the map $\varphi^{-1}:\varphi(U) \rightarrow U$ induces an algebra homomorphism $(\varphi^{-1})^{\ast}:\mathcal{F}_{p} \rightarrow \mathcal{F}_{\varphi(p)}$ given by $(\varphi^{-1})^{\ast}([f])=[f \circ \varphi^{-1}]$. Thus, the first order terms vanish, which means that for each $j=1,\ldots,m$
\begin{equation*}
\frac{\partial}{\partial t^{j}}\left ( \sum_{i=1}^{m}a_{i}t^{i} \right ) \biggl |_{0}=0
\end{equation*}
and therefore $a_{i}=0$, for all $i=1,\ldots,m$.\label{tangent_space}
\end{proof}
\end{prop}
From this we conclude that $V_{p}$ is finite dimensional and $dim(V_{p})=m$. An element $\xi_{p} \in V_{p}$ is called a tangent vector on $\manifold{M}$ at $p$.

For each $p \in \manifold{M}$ we can associate to each tangent vector $\xi_{i} \in V_{p}$ a linear map $v_{p}: \mathcal{F}_{p} \rightarrow \R$, given by
\begin{equation*}
v_{p}(f)=
\begin{cases}
	0 & \text{ , if } \exists \; c \in f \;| \; c(x)=c\;\forall x \in \manifold{M}\\
	\xi_{p}([f]) & \text{ , if } f \in I_{p}
\end{cases}
\end{equation*}
where $[f]$ denotes the class corresponding to the germ $f$ in $I_{p}/I_{p}^{2}$. Since all germs can be written as $f=\tilde{f}+f(p)$, where $\tilde{f} \in I_{p}$ and $f(p)$ is the germ of the constant function whose value is $f(p)$, $v_{p}$ satisfies the following property
\begin{eqnarray*}
v_{p}(fg) & = & v_{p}((\tilde{f}+f(p))(\tilde{g}+g(p)))=\\
& = & v_{p}(\tilde{f}\tilde{g}+f(p)\tilde{g}+g(p)\tilde{f}+f(p)g(p))=\\
& = & v_{p}(\tilde{f}\tilde{g})+v_{p}(f(p)\tilde{g})+v_{p}(g(p)\tilde{f})+v_{p}(f(p)g(p))=\\
& = & \xi_{p}(\tilde{f}\tilde{g})+f(p)\xi_{p}(\tilde{g})+g(p)\xi_{p}(\tilde{f})+0=\\
& = & f(p)\xi_{p}(\tilde{g})+g(p)\xi_{p}(\tilde{f})=\\
& = & f(p)v_{p}(\tilde{g})+g(p)v_{p}(\tilde{f})=\\
& = & g(p)v_{p}(f)+f(p)v_{p}(g)
\end{eqnarray*}
When a linear map $w:\mathcal{F}_{p} \rightarrow \R$ obeys such property we call $w$ a \emph{derivation on} $\mathcal{F}_{p}$ \emph{at the point} $p$.

Conversely, if $w$ is a derivation on $\mathcal{F}_{p}$ at $p$, we can associate $w$ to a unique tangent vector $\eta_{p}$ such that $\eta_{p}([f])=w(f)$ for all $f \in \mathcal{F}_{p}$. To see this, note that if $c$ is the germ of a constant function
\begin{equation*}
w(c)=w(c \cdot 1)=c w(1)=c w(1 \cdot 1)=c w(1)+c w(1)=2 w(c)
\end{equation*}
and therefore $w(c)=0$. By writing $f$ as $f=\tilde{f}+f(p)$, we have:
\begin{equation*}
w(f)=w(\tilde{f}+f(p))=w(\tilde{f})+w(f(p))=w(\tilde{f})
\end{equation*}
therefore the value of $w$ is determined by its value at $I_{p}$. If $f \in I_{p}^{2}$, there exists $g,h \in I_{p}$ such that $f=gh$. Hence
\begin{equation*}
w(f)=w(gh)=h(p)w(g)+g(p)w(h)=0
\end{equation*}
thus, $w$ vanishes on $I_{p}^{2}$. However, $f=\tilde{f}+f(p)$ gives
\begin{equation*}
w(f)=w(f-f(p))=w(\tilde{f})
\end{equation*}
which shows that if $\tilde{f} \; \equiv \; \tilde{g} \; mod \; I_{p}^{2}$, then $w(f)=w(g)$ and $w$ induces an unique linear map $\eta_{p}$ taking elements in $I_{p}/I_{p}^{2}$ to real values. In other words $\eta_{p} \in V_{p}$.

So we stablish an one to one correspondence between derivations on $\mathcal{F}_{p}$ at $p$ and elements in $(I_{p}/I_{p}^{2})^{\ast}$. It is not hard to see that the set of such derivations at a point, with the usual operations of addition and product by scalars, turns to be a $\R$-vector space and the association that sends elements in $(I_{p}/I_{p}^{2})^{\ast}$ to derivations on $\mathcal{F}_{p}$ at $p$ above mentioned is a $\R$-vector space isomorphism. Thus, we can speak in elements in $V_{p}$ acting on a germ $f \in \mathcal{F}_{p}$, once we implicitly understood the association above constructed. And beyond. Define the action of a tangent vector $v_{p} \in V_{p}$ on a function $f \in \funcman{\manifold{M}}$ by
\begin{equation*}
v_{p}(f)=v_{p}(\mathbf{f})
\end{equation*}
where $\mathbf{f}$ is the class of $f$ in $\mathcal{F}_{p}$. So, $v_{p}(g)=v_{p}(f)$ when $g \in \mathbf{f}$. Linearity and Leibniz rule follows straightforward from this definition.

From those facts, if $\manifold{M}$ is a $m$-dimensional differentiable manifold, we can construct the differentiable vector bundle $(\mathcal{E},\manifold{M},\pi,GL(\R^{m}))$ with total space $\mathcal{E}=\cup_{p \in \manifold{M}}V_{p}$, base space $\manifold{M}$, projection $\pi$ given by $\pi(v_{p})=p$ and typical fibre $\R^{m}$, with differentiable structure obtained from the structure on $\manifold{M}$. Indeed, let $(U,\varphi)$ be a local chart on the $m$-dimensional differentiable manifold $\manifold{M}$ around $p \in \manifold{M}$. By using the same notations of proposition \ref{tangent_space}, given $f \in \funcman{\manifold{M}}$, take elements $\frac{\partial}{\partial x^{i}} |_{p} \in V_{p}$ for $i=1,\ldots,m$, such that
\begin{equation*}
\frac{\partial}{\partial x^{i}} \biggl |_{p}(f)=\frac{\partial (f \circ \varphi^{-1})}{\partial t^{i}} \biggl |_{\varphi(p)}
\end{equation*}
Now, given $\eta \in \pi^{-1}(U)$, for all $f \in \funcman{\manifold{M}}$
\begin{equation*}
\eta(f) = \eta \left ( \sum_{i=1}^{m} \frac{\partial (\tilde{f} \circ \varphi^{-1})}{\partial t^{i}} \biggl |_{\varphi(\pi(\eta))} \mathbf{x}^{i} \right ) = \sum_{i=1}^{m}\frac{\partial f}{\partial x^{i}} \biggl |_{\pi(\eta)}\eta(x^{i})
\end{equation*}
allowing to write
\begin{equation}
\eta = \sum_{i=1}^{m} \eta(x^{i}) \frac{\partial}{\partial x^{i}} \biggl |_{\pi(\eta)}
\label{vector_coordinates}
\end{equation}
and we call this formula \emph{coordinate expression} of $\eta$ with respect to the local chart $(U,\varphi)$ and the values $\eta(x^{i})$ \emph{coordinates of} $\eta$. Hence we can define a map $\tilde{\varphi}:\pi^{-1}(U) \rightarrow \R^{m}$ given by
\begin{equation*}
\tilde{\varphi}(\eta)=(\eta(x^{1}),\ldots,\eta(x^{m}))
\end{equation*}
Finally, define the map $\phi:\pi^{-1}(U) \rightarrow \R^{m} \times \R^{m}$ given by
\begin{equation}
\phi(\eta)=((\varphi \circ \pi)(\eta),\tilde{\varphi}(\eta))
\end{equation}
Let $\atlas{A}{M}$ be the atlas of $\manifold{M}$. For each local chart $(U_{\alpha},\varphi_{\alpha})$ associate the map $\phi_{\alpha}:\pi^{-1}(U_{\alpha}) \rightarrow \R^{2m}$, constructed as above. Declare a set $V \subset \mathcal{E}$ open on $\mathcal{E}$ if and only if there exists an open set $V_{0}$ on $\R^{m}$ and an index $\alpha$ such that $V=\phi_{\alpha}^{-1}(V_{0})$. The collection of those sets is a base for the topology on $\mathcal{E}$ which makes $\mathcal{E}$ a topological manifold. Besides, if $(U_{\alpha},\varphi_{\alpha}),(U_{\beta},\varphi_{\beta}) \in \atlas{A}{M}$, then the map $\phi_{\beta} \circ \phi_{\alpha}^{-1}:\phi_{\alpha}(\pi^{-1}(U_{\alpha})) \rightarrow \phi_{\beta}(\pi^{-1}(U_{\beta}))$ is $C^{\infty}$. This shows that the collection $(\pi^{-1}(U_{\alpha}),\phi_{\alpha})$ defines a differentiable atlas on $\mathcal{E}$.

With this structures, we see that $(\mathcal{E},\manifold{M},\pi,GL(\R^{m}))$ is a differentiable vector bundle with typical fibre $\R^{m}$, total space $\mathcal{E}$, which is a $2m$-dimensional differentiable manifold, base space $\manifold{M}$, projection $\pi:\mathcal{E} \rightarrow \manifold{M}$, which is a differentiable surjective submersion, trivializations given by the maps $\phi_{\alpha}=(\varphi_{\alpha} \circ \pi,\tilde{\varphi}_{\alpha})$, structure group $GL(\R^{m})$ in which compatibility conditions of trivializations are satisfied, since map such as $\phi_{\beta} \circ \phi_{\alpha}^{-1}$ are diffeomorphisms. To fix notations, we will write $T\manifold{M}=\mathcal{E}$ and call $T\manifold{M}$ the \emph{tangent bundle} of $\manifold{M}$, since its elements can be faced as tangent vectors at points in $\manifold{M}$. From now on we will denote the vector space tangent at $p$ by $T_{p}\manifold{M}=V_{p}$.

We will use this
\begin{dfn}[Vector fields]
Let $\manifold{M}$ be a $m$-dimensional differentiable manifold. The $C^{\infty}$ sections of $\pi$ of $(T\manifold{M},\manifold{M},\pi,GL(\R^{m}))$ are called vector fields on $\manifold{M}$. We denote the $\R$-vector space of vector fields with the usual operations of addition and product by scalars pointwise by $\mathfrak{X}(\manifold{M})$. Also, $\mathfrak{X}$ has a $\funcman{\manifold{M}}$-module structure given by pointwise product by functions.\label{vector_fields}
\end{dfn}
This leads to the following
\begin{thm}
Let $\manifold{M}$ be a $m$-dimensional differentiable manifold and $\funcman{\manifold{M}}$ the $\R$-algebra of $C^{\infty}$ functions on $\manifold{M}$. Then
\begin{equation*}
\mathfrak{X}(\manifold{M}) \approx_{\catvect{R}} Der(\funcman{\manifold{M}})
\end{equation*}\label{iso_fields_derivations}
\end{thm}
The reader is invited to proof the above theorem using the tools (and definitions) given here as an exercise.

Following this steps, we can now construct a differentiable vector bundle over a differentiable manifold $\manifold{M}$, whose differentiable sections corresponds to derivations of order $\leq$ to $r$ on the algebra of functions $\funcman{\manifold{M}}$ (definition \ref{high_order_derivation}).

First, we will need the notion of high order derivation at a point. Let $p \in \manifold{M}$ and $\mathcal{F}_{p}$ be the $\R$-algebra of germs of functions at $p$. If $I_{p}$ denotes the ideal of germs of functions vanishing at $p$ we have (as above) a natural $\R$-vector space structure on $I_{p}/I_{p}^{r}$, since $I_{p}^{r}$ (here $r \in \mathbb{Z},\;r \geq 1$) is an ideal of $I_{p}$. By denoting $J_{p}^{r}=(I_{p}/I_{p}^{r+1})^{\ast}$, we can repeat all what we have done on $V_{p}$ and define derivation of order $\leq$ $r$ at $p$.
\begin{dfn}
Let $\mathcal{F}_{p}$ be the germs of functions at a point $p$ in a $m$-dimensional differentiable manifold $\manifold{M}$. We say that a $\R$-linear map $D_{p}:\mathcal{F}_{p} \rightarrow \R$ is a differential operator of order $\leq$ $r$, $r \geq 1$, at $p$ if and only if for all $g \in \mathcal{F}_{p}$, the map $d_{g}:\mathcal{F}_{p} \rightarrow \R$ given by
\begin{displaymath}
d_{g}(f)=D_{p}(gf)-g(p)D_{p}(f)
\end{displaymath}
is a differential operator of order $\leq$ $r-1$ at $p$, and called a differential operator of order 0 at $p$ if it is a product by a germ of functions at $p$. $D_{p}$ is called a derivation of order $\leq$ $r$ at $p$ if, in addiction, it is identically null on germs of constant functions.
\end{dfn}
Compare with the definition \ref{high_order_derivation}.

$J_{p}^{r}$ is finite dimensional. To see this, let $(U,\varphi)$ be a local chart around $p \in \manifold{M}$. Without loss of generality, we can suppose $\varphi(U)$ convex and $\varphi(p)=0$. As before, denote the coordinates on $U$ by $x^{i}=t^{i} \circ \varphi$, where $t^{i}:\R^{m}\rightarrow \R$ is the $i$-th projection on $\R^{m}$. The claim follows by showing that the set of classes of the functions $x^{i},\;i=1,\ldots,m$, $x^{i}x^{j},1 \leq i \leq j \leq n,\;\ldots, x^{i_{1}} \cdot \ldots \cdot x^{i_{r}},\;i_{1} \leq \ldots \leq i_{r}$ is a basis for $I_{p}/I_{p}^{r+1}$. Let $f \in I_{p}/I_{p}^{r+1}$. Let $f \in \funcman{\manifold{M}}$ be a representing element of this class. By the Taylor formula, the coordinate expression of $f$ on $U$ is written by
\begin{eqnarray*}
f & = & \sum_{i=1}^{m}\frac{\partial (f \circ \varphi^{-1})}{\partial t^{i}}\biggl |_{0}x^{i} + \frac{1}{2}\sum_{i,j=1}^{m}\frac{\partial^{2}(f \circ \varphi^{-1})}{\partial t^{i} \partial t^{j}}\biggl |_{0}x^{i}x^{j} + \ldots + \\
& + & \frac{1}{r!}\sum_{i_{1},\ldots,i_{r+1}=1}^{m}x^{i_{1}}\ldots x^{i_{r+1}}\int_{0}^{1}(1-s)^{r}\frac{\partial^{r+1}(f \circ \varphi^{-1})}{\partial t^{i_{1}}\ldots \partial t^{i_{r+1}}}\biggl |_{sy}ds
\end{eqnarray*}
By taking the quotient, we note that the last term on the right hand side vanishes on $I_{p}/I_{p}^{r+1}$. Hence, the class of $f$ is written
\begin{eqnarray*}
f & = & \sum_{i=1}^{m}\frac{\partial (f \circ \varphi^{-1})}{\partial t^{i}}\biggl |_{0}x^{i} + \frac{1}{2}\sum_{i,j=1}^{m}\frac{\partial^{2}(f \circ \varphi^{-1})}{\partial t^{i} \partial t^{j}}\biggl |_{0}x^{i}x^{j} + \ldots + \\
& + & \frac{1}{r!}\sum_{i_{1},\ldots,i_{r}=1}^{m} \frac{\partial^{r}(f \circ \varphi^{-1})}{\partial t^{i_{1}}\ldots \partial t^{i_{r}}}\biggl |_{0} x^{i_{1}}\ldots x^{i_{r}}
\end{eqnarray*}
which shows that $x^{i},\;i=1,\ldots,m$, $x^{i}x^{j},1 \leq i \leq j \leq n,\;\ldots, x^{i_{1}} \cdot \ldots \cdot x^{i_{r}},\;i_{1} \leq \ldots \leq i_{r}$ spans $I_{p}/I_{p}^{r+1}$, since $f\circ \varphi^{-1}$ is differentiable.

To show the linear independence, note that
\begin{eqnarray*}
& & \sum_{i=1}^{m}a_{i}x^{i} + \sum_{1 \leq i \leq j \leq m}a_{ij}x^{i}x^{j} + \ldots + \sum_{i_{1} \leq \ldots \leq i_{r}}a_{i_{1}\ldots i_{r}}x^{i_{1}}\ldots x^{i_{r}} = 0 \Rightarrow \\
& & \Rightarrow \sum_{i=1}^{m}a_{i}x^{i} + \sum_{1 \leq i \leq j \leq m}a_{ij}x^{i}x^{j} + \ldots + \sum_{i_{1} \leq \ldots \leq i_{r}}a_{i_{1}\ldots i_{r}}x^{i_{1}}\ldots x^{i_{r}} \in I_{p}^{r+1} \Rightarrow \\
& & \Rightarrow \sum_{i=1}^{m}a_{i}t^{i} + \sum_{1 \leq i \leq j \leq m}a_{ij}t^{i}t^{j} + \ldots + \sum_{i_{1} \leq \ldots \leq i_{r}}a_{i_{1}\ldots i_{r}}t^{i_{1}}\ldots t^{i_{r}} \in I_{\varphi(p)}^{r+1}
\end{eqnarray*}
Thus, the terms of order $\leq$ $r$ are null, which leads to
\begin{eqnarray*}
& & \frac{\partial^{r}}{\partial t^{j_{1}} \ldots \partial t^{j_{r}}} \left ( \sum_{i_{1} \leq \ldots \leq i_{r}}a_{i_{1}\ldots i_{r}}t^{i_{1}}\ldots t^{i_{r}} \right )\biggl |_{0} = 0\\
& & \vdots \\
& & \frac{\partial^{2}}{\partial t^{k} \partial t^{l}} \left ( \sum_{1 \leq i \leq j \leq m}a_{ij}t^{i}t^{j} \right ) \biggl |_{0} = 0 \\
& & \frac{\partial}{\partial t^{j}} \left ( \sum_{i=1}^{m}a_{i}t^{i} \right ) \biggl |_{0} = 0
\end{eqnarray*}
hence
\begin{eqnarray*}
& & a_{j_{1}\ldots j_{r}} = 0 \\
& & \vdots \\
& & a_{kl} = 0 \\
& & a_{j} = 0
\end{eqnarray*}
for all possible combinations of indices. So, $I_{p}/I_{p}^{r+1}$ is finite dimensional, therefore $J_{p}^{r}$ is finite dimensional also.

Let $\xi_{p} \in J_{p}^{r}$. Associate to $\xi_{p}$ a linear map $D_{p}:\mathcal{F}_{p} \rightarrow \R$ given by
\begin{equation*}
D_{p}(f)=
\begin{cases}
	0 & \text{ , if } \exists \; c \in f \;| \; c(x)=c\;\forall x \in \manifold{M} \\
	\xi_{p}([f]) & \text{ , if } f \in I_{p}
\end{cases}
\end{equation*}
where $[f]$ denotes the class of $f$ in $I_{p}/I_{p}^{r+1}$. By writing germs $f \in \mathcal{F}_{p}$ as $f=\tilde{f}+f(p)$, with $\tilde{f} \in I_{p}$, given $g \in \mathcal{F}_{p}$ we have
\begin{eqnarray}
& & \Delta_{g}(f) = D_{p}(gf)-g(p)D_{p}(f) = \nonumber \\
& & = D_{p}(\tilde{g}\tilde{f})+g(p)D_{p}(\tilde{f})+f(p)D_{p}(\tilde{g})+2f(p)g(p)D_{p}(1)-g(p)D_{p}(\tilde{f}) = \nonumber \\
& & = D_{p}(\tilde{g}\tilde{f})+f(p)D_{p}(\tilde{g}) = \xi_{p}(\tilde{g}\tilde{f})+f(p)\xi_{p}(\tilde{g}).
\label{proof_jets_order_r_op_at_p}
\end{eqnarray}
Note that, given $f_{1} \in I_{p}$, for all $f_{0} \in I_{p}$, we have
\begin{displaymath}
\delta^{r-1}_{f_{1}}(f_{0})=\xi_{p}(f_{1}f_{0})-f_{1}(p)\xi_{p}(f_{0})=\xi_{p}(f_{1}f_{0})
\end{displaymath}
and successively we can see that, given $f_{1},\ldots,f_{i}$
\begin{displaymath}
\delta^{r-i}_{f_{i}}(f_{0})=\xi_{p}(f_{i}f_{i-1} \ldots f_{0})
\end{displaymath}
for all $f_{0} \in I_{p}$. Now
\begin{displaymath}
\delta^{0}_{f_{r}}(f_{0})=\xi_{p}(f_{r} \ldots f_{0}) = 0
\end{displaymath}
which shows that $\delta^{1}_{f_{r-1}}$ is a differential operator of order $\leq$ 1 at the point $p$, viewed as restricted to $I_{p}$. Restricting to $I_{p}$, by construction, $\delta^{r-i}_{f_{i}}$ being a differential operator of order $\leq$ $r-i$ at $p$ implies that $\delta^{r-i+1}_{f_{i-1}}$ is a differential operator of order $\leq$ to $r-i+1$ at $p$. Hence, we have $\xi_{p}$ differential operator of order $\leq$ to $r$ at $p$ and also, if $\delta:I_{p} \rightarrow \R$ is an operator such that for all $f \in I_{p}$, $\delta(f)=\xi_{p}(f_{1} \ldots f_{k}f)$, with $f_{1}, \ldots, f_{k} \in I_{p}$, then $\delta$ is a differential operator of order $\leq$ $r-k$ at $p$. Putting on equation \ref{proof_jets_order_r_op_at_p}, we see that $\Delta_{g}$ is a differential operator of order $\leq$ $r-1$ at $p$, leading to the conclusion that $D_{p}$ is a differential operator of order $\leq$ $r$ at the point $p$.

Conversely, let $\omega:\mathcal{F}_{p} \rightarrow \R$ be a derivation of order $\leq$ $r$ at the point $p$. Then $f \in \mathcal{F}_{p}$ gives
\begin{displaymath}
\omega(f)=\omega(\tilde{f}+f(p))=\omega(\tilde{f})
\end{displaymath}
which shows that the value of $\omega$ depends on its evaluation at $I_{p}$ only. We also have that if $f \in I_{p}^{r+1}$ then there exists $f_{1},\ldots,f_{r+1} \in I_{p}$ such that $f=f_{1}\ldots f_{r+1}$ and therefore
\begin{eqnarray*}
& & \omega(f)=\omega(f_{1}\ldots f_{r+1})=\delta^{r-1}_{f_{r+1}}(f_{1}\ldots f_{r})+f_{r+1}(p)\omega(f_{1}\ldots f_{r}) = \\
& & = \delta^{r-1}_{f_{r+1}}(f_{1}\ldots f_{r}) = \delta^{r-2}_{f_{r}}(f_{1}\ldots f_{r-1}) = \ldots = \delta^{0}_{f_{2}}(f_{1}) = \\
& & = f_{1}(p)g = 0
\end{eqnarray*}
for some $g \in \mathcal{F}_{p}$, where $\delta^{r-i}_{f_{r-i+2}}$ is a differential operator of order $\leq$ $r-i$, for $i=1,\ldots,r$. So, if $f \equiv g$ $mod\;I_{p}^{r+1}$ in $I_{p}$ then $\omega(f)=\omega(g)$ and $\omega$ can be viewed as an element in $(I_{p}/I_{p}^{r+1})^{\ast}$. It follows that there exists an one to one correspondence between derivations of order $\leq$ $r$ at the point $p$ and elements in $J_{p}^{r}$.

We define the action of a derivation of order $\leq$ $r$ at the point $p$, $D_{p}$, on a function $f \in \funcman{\manifold{M}}$ as given by
\begin{displaymath}
D_{p}(f)=\xi_{p}([f])
\end{displaymath}
where $\xi_{p}$ is the element related to $D_{p}$ by the above correspondence and $[f]$ is the equivalence class of the function $f$ in $I_{p}/I_{p}^{r+1}$.

We can now prove the following theorem.
\begin{thm}
Let $\manifold{M}$ be a $m$-dimensional differentiable manifold. There exists a differentiable vector bundle $J^{r}(\manifold{M})$, whose base space is $\manifold{M}$ and whose space of differentiable sections $\Gamma(J^{r}(\manifold{M}))$ is isomorphic, as $\R$-vector space, to the space of derivations of order $\leq$ $r$ on $\funcman{\manifold{M}}$.
\begin{proof}
Let us take the differentiable vector bundle $\displaystyle J^{r}(\manifold{M})=\bigcup_{p \in \manifold{M}}J_{p}^{r}$, with $J_{p}^{r}=(I_{p}/I_{p}^{r+1})^{\ast}$ whose coordinate functions $\phi:\pi^{-1}(U) \rightarrow \R^{K}$, with $\pi$ the projection of $J^{r}(\manifold{M})$ on $\manifold{M}$ and $U$ an open subset of $\manifold{M}$, are of the form $\phi(\omega_{p}) = (x^{i}(\pi(\omega_{p})),\omega_{p}(x^{i}),\omega_{p}(x^{i}x^{j}),\ldots,\omega_{p}(x^{i_{1}}\ldots x^{i_{r}}))$, where the indices are increasing, 1 to $m$, for all $\omega_{p} \in J_{p}^{r}$. Note that $K=\sum_{k=1}^{r}\binom{m+k-1}{k}$.

Let $\omega:\manifold{M}\rightarrow J^{r}(\manifold{M})$ a differentiable section of $J^{r}(\manifold{M})$. We associate $\omega$ to the map $D:\funcman{\manifold{M}}\rightarrow \funcman{\manifold{M}}$ given by
\begin{displaymath}
D(f)(p)=\omega(p)(f)\quad \forall f \in \funcman{\manifold{M}}
\end{displaymath}

The map $D$ above defined is a derivation of order $\leq$ $r$ on $\funcman{\manifold{M}}$. To see this, first note that if $c$ is a constant function, then
\begin{displaymath}
D(c)(p)=\omega(p)(c)=0\quad \forall p \in \manifold{M}
\end{displaymath}
hence $D$ vanishes on constants. Second, given $g \in \funcman{\manifold{M}}$, the operator $\Delta_{g}$ given by
\begin{displaymath}
\Delta_{g}(f)=D(gf)-gD(f)\quad \forall f \in \funcman{\manifold{M}}
\end{displaymath}
is such that
\begin{eqnarray*}
\Delta_{g}(f)(p) & = & D(gf)(p)-g(p)D(f)(p) = \omega(p)(gf)-g(p)\omega(p)(f) = \\
& = & \delta_{g}(p)(f)
\end{eqnarray*}
which is a differentiable section (by construction) of $J^{r-1}(\manifold{M})$. However, $\Gamma(J^{1}(\manifold{M}))=\mathfrak{X}(\manifold{M})$ and theorem \ref{iso_fields_derivations} shows that $\Gamma(J^{1}(\manifold{M})) \approx_{\catvect{\R}} Der(\funcman{\manifold{M}})$. Therefore, by induction, $D$ is a derivation of order $\leq$ $r$ on $\funcman{\manifold{M}}$.

The assignment $\omega \mapsto D$ is clearly linear. Lets show it is injective. Suppose $\omega$ associated to $D$ identically null. We have
\begin{eqnarray*}
D(f) = 0 & , & \forall f \in \funcman{\manifold{M}} \\
D(f)(p) = 0 & , & \forall f \in \funcman{\manifold{M}},\;\forall p \in \manifold{M} \\
\omega(p)(f) = 0 & , & \forall p \in \manifold{M},\;\forall f \in \funcman{\manifold{M}} \\
\omega(p) = 0 & , & \forall p \in \manifold{M} \\
\omega = 0
\end{eqnarray*}
Lets show it is surjective. Let $D:\funcman{\manifold{M}} \rightarrow \funcman{\manifold{M}}$ be a derivation of order $\leq$ $r$ on $\funcman{\manifold{M}}$. For each $p \in \manifold{M}$, define $\omega_{p}:\mathcal{F}_{p} \rightarrow \R$ by
\begin{displaymath}
\omega_{p}(f) = D(f)(p) \quad \forall f \in \funcman{\manifold{M}}
\end{displaymath}
where $f$ on the left hand side is for the equivalence class in $\mathcal{F}_{p}$ of the function represented by the symbol $f$ on the right hand side. $\omega_{p}$ is well defined because if $f,g \in \funcman{\manifold{M}}$ are such that $f \equiv g$ in $\mathcal{F}_{p}$, we can take a local chart $(U,\varphi)$ around $p$ such that $U \subset W$, where $W$ is an open subset of $\manifold{M}$ in which $f$ and $g$ coincides, $\varphi(p)=0$ and $\varphi(U)$ is convex. Now, on $U$, $f$ and $g$ are written
\begin{eqnarray*}
f & = & f(p) + \sum_{i=1}^{m}\frac{\partial (f \circ \varphi^{-1})}{\partial t^{i}}\biggl |_{0}x^{i} + \frac{1}{2}\sum_{i,j=1}^{m}\frac{\partial^{2}(f \circ \varphi^{-1})}{\partial t^{i} \partial t^{j}}\biggl |_{0}x^{i}x^{j} + \ldots + \\
& + & \frac{1}{r!}\sum_{i_{1},\ldots,i_{r+1}=1}^{m}x^{i_{1}}\ldots x^{i_{r+1}}\int_{0}^{1}(1-s)^{r}\frac{\partial^{r+1}(f \circ \varphi^{-1})}{\partial t^{i_{1}}\ldots \partial t^{i_{r+1}}}\biggl |_{sy}ds
\end{eqnarray*}
and
\begin{eqnarray*}
g & = & g(p) + \sum_{i=1}^{m}\frac{\partial (g \circ \varphi^{-1})}{\partial t^{i}}\biggl |_{0}x^{i} + \frac{1}{2}\sum_{i,j=1}^{m}\frac{\partial^{2}(g \circ \varphi^{-1})}{\partial t^{i} \partial t^{j}}\biggl |_{0}x^{i}x^{j} + \ldots + \\
& + & \frac{1}{r!}\sum_{i_{1},\ldots,i_{r+1}=1}^{m}x^{i_{1}}\ldots x^{i_{r+1}}\int_{0}^{1}(1-s)^{r}\frac{\partial^{r+1}(g \circ \varphi^{-1})}{\partial t^{i_{1}}\ldots \partial t^{i_{r+1}}}\biggl |_{sy}ds
\end{eqnarray*}
where $t^{i}:\R^{m}\rightarrow \R$ is the $i$-th canonical projection on $\R^{m}$. Hence, on $U$,
\begin{eqnarray*}
D(f) & = & \sum_{i=1}^{m}\frac{\partial (f \circ \varphi^{-1})}{\partial t^{i}}\biggl |_{0}D(x^{i}) + \frac{1}{2}\sum_{i,j=1}^{m}\frac{\partial^{2}(f \circ \varphi^{-1})}{\partial t^{i} \partial t^{j}}\biggl |_{0}D(x^{i}x^{j}) + \ldots + \\
& + & \frac{1}{r!}\sum_{i_{1},\ldots,i_{r+1}=1}^{m}D(x^{i_{1}}\ldots x^{i_{r+1}})\int_{0}^{1}(1-s)^{r}\frac{\partial^{r+1}(f \circ \varphi^{-1})}{\partial t^{i_{1}}\ldots \partial t^{i_{r+1}}}\biggl |_{sy}ds
\end{eqnarray*}
and
\begin{eqnarray*}
D(g) & = & \sum_{i=1}^{m}\frac{\partial (g \circ \varphi^{-1})}{\partial t^{i}}\biggl |_{0}D(x^{i}) + \frac{1}{2}\sum_{i,j=1}^{m}\frac{\partial^{2}(g \circ \varphi^{-1})}{\partial t^{i} \partial t^{j}}\biggl |_{0}D(x^{i}x^{j}) + \ldots + \\
& + & \frac{1}{r!}\sum_{i_{1},\ldots,i_{r+1}=1}^{m}D(x^{i_{1}}\ldots x^{i_{r+1}})\int_{0}^{1}(1-s)^{r}\frac{\partial^{r+1}(g \circ \varphi^{-1})}{\partial t^{i_{1}}\ldots \partial t^{i_{r+1}}}\biggl |_{sy}ds
\end{eqnarray*}
But since $D$ is a derivation of order $\leq$ $r$, we have
\begin{displaymath}
D(x^{i_{1}}\ldots x^{i_{r+1}})(p) = 0
\end{displaymath}
for all relevant combination of indices. As $f$ and $g$ are in the same germ of functions at $p$, all partial derivatives up to order $r$ of its coordinate expressions coincide, leading to
\begin{eqnarray*}
& & \omega_{p}(f) = D(f)(p) = \\
& & = \sum_{i=1}^{m}\frac{\partial (f \circ \varphi^{-1})}{\partial t^{i}}\biggl |_{0}D(x^{i})(p) + \frac{1}{2}\sum_{i,j=1}^{m}\frac{\partial^{2}(f \circ \varphi^{-1})}{\partial t^{i} \partial t^{j}}\biggl |_{0}D(x^{i}x^{j})(p) + \ldots + \\
& & + \frac{1}{r!}\sum_{i_{1},\ldots,i_{r+1}=1}^{m}D(x^{i_{1}}\ldots x^{i_{r+1}})(p)\int_{0}^{1}(1-s)^{r}\frac{\partial^{r+1}(f \circ \varphi^{-1})}{\partial t^{i_{1}}\ldots \partial t^{i_{r+1}}}\biggl |_{sy}ds = \\
& & + \sum_{i=1}^{m}\frac{\partial (g \circ \varphi^{-1})}{\partial t^{i}}\biggl |_{0}D(x^{i})(p) + \frac{1}{2}\sum_{i,j=1}^{m}\frac{\partial^{2}(g \circ \varphi^{-1})}{\partial t^{i} \partial t^{j}}\biggl |_{0}D(x^{i}x^{j})(p) + \ldots + \\
& + & \frac{1}{r!}\sum_{i_{1},\ldots,i_{r+1}=1}^{m}D(x^{i_{1}}\ldots x^{i_{r+1}})(p)\int_{0}^{1}(1-s)^{r}\frac{\partial^{r+1}(g \circ \varphi^{-1})}{\partial t^{i_{1}}\ldots \partial t^{i_{r+1}}}\biggl |_{sy}ds = \\
& & = D(g)(p) = \omega_{p}(g)
\end{eqnarray*}
As $D$ is a derivation of order $\leq$ $r$, it follows by induction on $r$ that $\omega_{p}$ is a derivation of order $r$ at $p$. Hence, $\omega_{p} \in J_{p}^{r}$, for all $p \in \manifold{M}$. Let $\omega:\manifold{M}\rightarrow J^{r}(\manifold{M})$ be the map given by $\omega(p)=\omega_{p}$. For all $f \in \funcman{\manifold{M}}$, we have
\begin{displaymath}
\omega(p)(f) = \omega_{p}(f) = D(f)(p)
\end{displaymath}
showing that $p \mapsto \omega(p)(f)$ is differentiable, because $D(f)$ is, and so $\omega$ is a differentiable section of $J^{r}(\manifold{M})$.

Hence, $\Gamma(J^{r}(\manifold{M}))$ is isomorphic, as $\R$-vector space, to the space of derivations of order $\leq$ $r$ on $\funcman{\manifold{M}}$.\label{iso_deriv_order_r_jets}
\end{proof}
\end{thm}

It follows from the last theorem that if $\manifold{M}$ is a $m$-dimensional manifold, a derivation of order $\leq$ $r$, $D:\funcman{\manifold{M}}\rightarrow \funcman{\manifold{M}}$ is written locally as
\begin{displaymath}
D(f) = \sum_{k=1}^{r} \sum_{1 \leq i_{1} \leq \ldots \leq i_{k} \leq m}a_{i_{1}\ldots i_{k}}(x_{1},\ldots,x_{m})\frac{\partial^{k}f}{\partial x^{i_{1}}\ldots \partial x^{i_{k}}},
\end{displaymath}
with $a_{i_{1}\ldots i_{k}}$ differentiable \cite{kostrikin_safarevich_algebra_1}.

\section{Multidifferential Operators}
Let $(A,\mu,e)$ be an associative commutative unital $\field{K}$-algebra. Denote
\begin{displaymath}
C^{n}(A,A)=Hom_{\catvect{K}}(A^{\otimes n},A),\quad \forall n \in \mathbb{Z}, n \geq 0
\end{displaymath}
\begin{displaymath}
C^{\bullet}(A,A)=\bigoplus_{n \geq 0}C^{n}(A,A)
\end{displaymath}
\begin{dfn}[Partial composition]
Let $f \in C^{m+1}(A,A)$ and $g \in C^{n+1}(A,A)$. For $1 \leq i \leq m+1$, the $i$-th partial composition of $f$ and $g$ is the linear map $\circ_{i}:C^{m+1}(A,A) \otimes C^{n+1}(A,A) \rightarrow C^{m+n+1}(A,A)$ given by
\begin{displaymath}
f \circ_{i} g = f(id_{A}^{\otimes (i-1)} \otimes g \otimes id_{A}^{\otimes (m-i+1)})
\end{displaymath}
where $id_{A}$ denotes the identity of $A$.\label{partial_composition}
\end{dfn}
\begin{dfn}[Total composition]
The total composition $\circ :C^{\bullet}(A,A) \otimes C^{\bullet}(A,A) \rightarrow C^{\bullet}(A,A)$ is the linear map which, for each $f \in C^{m+1}(A,A)$ and $g \in C^{n+1}(A,A)$, associates $f \circ g \in C^{m+n+1}(A,A)$ given by
\begin{displaymath}
f \circ g = \sum_{i=1}^{m+1}(-1)^{n(i+1)}f \circ_{i} g.
\end{displaymath}\label{total_composition}
\end{dfn}
\begin{dfn}[Cup product]
The cup product is the degree 0 $\field{K}$-linear map $\smile:C^{\bullet}(A,A) \otimes C^{\bullet}(A,A) \rightarrow C^{\bullet}(A,A)$, which for each $f \in C^{m+1}(A,A)$ and $g \in C^{n+1}(A,A)$, associates
\begin{displaymath}
f \smile g = (-1)^{(m+1)(n+1)} \mu \circ (f \otimes g)
\end{displaymath}
\emph{i.e.}, if $a_{0},\ldots,a_{m},a_{m+1},\ldots,a_{m+n+1} \in A$, we have
\begin{eqnarray*}
& & (f \smile g)(a_{0} \otimes \ldots \otimes a_{m} \otimes a_{m+1} \otimes \ldots \otimes a_{m+n+1}) = \\
& & = \mu(f(a_{0} \otimes \ldots \otimes a_{m}) \otimes g(a_{m+1} \otimes \ldots \otimes a_{m+n+1})).\label{cup_product}
\end{eqnarray*}
\end{dfn}
\begin{prop}
$(C^{\bullet}(A,A),\smile)$ is an associative graded $\field{K}$-algebra.
\begin{proof}
By construction, $C^{\bullet}(A,A)$ is a graded $\field{K}$-vector space. By $\field{K}$-linearity and 0 degree of cup product, we only have to prove the associativity. Note that $\mu \in C^{2}(A,A)$, which leads to
\begin{displaymath}
\mu^{2} = \mu \circ \mu = \mu(\mu \otimes id_{A} - id_{A} \otimes \mu) = 0
\end{displaymath}
So, if $f \in C^{m+1}(A,A),\;g \in C^{n+1}(A,A),\;h \in C^{l+1}(A,A)$, and denoting $\sigma = (m+1)(n+1)+(m+1)(l+1)+(n+1)(l+1)$ we have
\begin{displaymath}
(f \smile g) \smile h - f \smile (g \smile h) = (-1)^{\sigma}\mu^{2}(f \otimes g \otimes h) = 0 \qedhere\label{cup_product_associativity}
\end{displaymath}
\end{proof}
\end{prop}
\begin{dfn}[Hochschild cohomology]
The Hochschild cohomology of an associative $\field{K}$-algebra $A$ with coefficients in $A$ is the cohomology of the complex
\begin{displaymath}
\xymatrix@1{0 \; \ar[r] & \; A \; \ar[r]^-{\delta_{H}} & \; C^{1}(A,A) \; \ar[r]^-{\delta_{H}} & \; \ldots \; \ar[r]^-{\delta_{H}} & \; C^{n}(A,A) \; \ar[r]^-{\delta_{H}} & \; \ldots }
\end{displaymath}
where the coboundary operator $\delta_{H}$, called the Hochschild differential, is given by
\begin{eqnarray*}
& & (\delta_{H}f)(a_{0} \otimes \ldots \otimes a_{n}) := \mu(a_{0} \otimes f(a_{1} \otimes \ldots \otimes a_{n})) + \\
& & + \sum_{i=0}^{n-1}(-1)^{i+1}f(a_{0} \otimes \ldots \otimes \mu(a_{i} \otimes a_{i+1}) \otimes \ldots \otimes a_{n}) + \\
& & + (-1)^{n+1}\mu(f(a_{0} \otimes \ldots \otimes a_{n-1}) \otimes a_{n})\end{eqnarray*}\label{hochschild_complex}
for any $f \in C^{n}(A,A)$, for all $a_{i} \in A,\; i=0,\ldots,n$.
\end{dfn}
\begin{prop}
The Hochschild differential $\delta_{H}$ is a degree 1 derivation for $(C^{\bullet}(A,A),\smile)$.
\begin{proof}
For simplicity, we will denote the product $\mu$ of the algebra $A$ by juxtaposition. By linearity, it is enough to consider the evaluation of $\delta_{H}$ on products of homogeneous terms. Let $f \in C^{m+1}(A,A)$ and $g \in C^{n+1}(A,A)$. For any $a_{i} \in A,\; i=0,\ldots,m+n+2$, we have
\begin{eqnarray*}
& & \delta_{H}(f \smile g)(a_{0} \otimes \ldots \otimes a_{m+n+2}) = a_{0}(f \smile g)(a_{1} \otimes \ldots \otimes a_{m+n+2}) + \\
& + & \sum_{i=0}^{m+n+1}(-1)^{i+1}(f \smile g)(a_{0} \otimes \ldots \otimes a_{i}a_{i+1} \otimes \ldots \otimes a_{m+n+2}) + \\
& + & (-1)^{m+n+3}(f \smile g)(a_{0} \otimes  \ldots  \otimes a_{m+n+1})a_{m+n+2} = \\
& = & a_{0}f(a_{1}\otimes \ldots \otimes a_{m+1})g(a_{m+2}\otimes \ldots \otimes a_{m+n+2}) + \\
& + & \sum_{i=0}^{m}(-1)^{i+1}f(a_{0}\otimes \ldots \otimes a_{i}a_{i+1}\otimes \ldots \otimes a_{m+1})g(a_{m+2}\otimes \ldots \otimes a_{m+n+2}) + \\
& + & \sum_{i=m+1}^{m+n+1}(-1)^{i+1}f(a_{0}\otimes \ldots \otimes a_{m})g(a_{m+1}\otimes \ldots \otimes a_{i}a_{i+1}\otimes \ldots \otimes a_{m+n+2}) + \\
& + & (-1)^{m+n+3}f(a_{0}\otimes \ldots \otimes a_{m})g(a_{m+1}\otimes \ldots \otimes a_{m+n+1})a_{m+n+2} = \\
& = & (a_{0}f(a_{1}\otimes \ldots \otimes a_{m}) + \sum_{i=0}^{m}(-1)^{i+1}f(a_{0}\otimes \ldots \otimes a_{i}a_{i+1}\otimes \ldots \otimes a_{m+1}) + \\
& + & (-1)^{m+2}f(a_{0}\otimes \ldots \otimes a_{m})a_{m+1})g(a_{m+2}\otimes \ldots \otimes a_{m+n+2}) + \\
& + & (-1)^{m+1}f(a_{0}\otimes \ldots \otimes a_{m})(a_{m+1}g(a_{m+2}\otimes \ldots \otimes a_{m+n+2}) + \\
& + & \sum_{i=0}^{n}(-1)^{i+1}g(a_{m+1} \otimes \ldots \otimes a_{m+i+1}a_{m+i+2} \otimes \ldots \otimes a_{m+n+2}) + \\
& + & (-1)^{n+2}g(a_{m+1}\otimes \ldots \otimes a_{m+n+1})a_{m+n+2}) = \\
& = & ((\delta_{H}f)\smile g)(a_{0}\otimes \ldots \otimes a_{m+n+2}) + (-1)^{m+1}(f \smile \delta_{H}g)(a_{0}\otimes \ldots \otimes a_{m+n+2})
\end{eqnarray*}
\end{proof}
\end{prop}
\begin{dfn}[Gerstenhaber bracket]
The Gerstenhaber bracket is the degree -1 $\field{K}$-linear map $\liebr{\;}{\;}:C^{\bullet}(A,A) \otimes C^{\bullet}(A,A) \rightarrow C^{\bullet}(A,A)$ which, for each $f \in C^{m+1}(A,A)$ and $g \in C^{n+1}(A,A)$
\begin{displaymath}
\liebr{f}{g}=f \circ g - (-1)^{mn}g \circ f.\label{gerstenhaber_bracket}
\end{displaymath}
\end{dfn}
\begin{prop}
$(C^{\bullet}(A,A),\liebr{\;}{\;})$ is a Lie superalgebra with respect to the reduced (by one) degree.
\end{prop}
The reader can find a proof of this fact in \cite{gerstenhaber_cohomology_structure}. 
\begin{prop}
Let $(A,\mu)$ be an associative $\field{K}$-algebra. For any $f \in C^{m+1}(A,A)$
\begin{displaymath}
\delta_{H}(f) = (-1)^{m}\liebr{\mu}{f}
\end{displaymath}
where $\liebr{\;}{\;}$ is the Gerstenhaber bracket.
\begin{proof}
Let $f \in C^{m+1}(A,A)$. Since $\mu \in C^{2}(A,A)$, it follows that
\begin{eqnarray*}
\liebr{\mu}{f} & = & \mu \circ f -(-1)^{m}f \circ \mu = \mu(f \otimes id_{A}) + (-1)^{m}\mu(id_{A} \otimes f) + \\
& + & (-1)^{m}\sum_{i=0}^{m}(-1)^{i+1}f(id_{A}^{\otimes i} \otimes \mu \otimes id_{A}^{\otimes (m-i)}) = \\
& = & (-1)^{m}\delta_{H}(f) \qedhere
\end{eqnarray*}
\end{proof}\label{gerstenhaber_bracket_hochschild_diff}
\end{prop}
\begin{prop}
Let $A$ be a $\field{K}$-algebra, where $\field{K}$ has characteristic different from 2. Fix a product $\nu \in C^{2}(A,A)$. Then $\nu$ is associative if and only if $\liebr{\nu}{\nu}=0$.
\begin{proof}
Given $f \in C^{m+1}(A,A)$, we have
\begin{eqnarray*}
\delta_{H}^{2}(f) & = & \delta_{H}(\delta_{H}(f)) = \\
& = & \delta_{H}((-1)^{m}\liebr{\nu}{f}) = (-1)^{2m+2}\liebr{\nu}{\liebr{\nu}{f}} = \liebr{\liebr{\nu}{\nu}}{f}-\liebr{\nu}{\liebr{\nu}{f}} \therefore \\
& & \delta_{H}^{2}(f) = \frac{1}{2}\liebr{\liebr{\nu}{\nu}}{f}
\end{eqnarray*}
To finish the proof, we just have to note that
\begin{displaymath}
\frac{1}{2}\liebr{\nu}{\nu} = \nu(\nu \otimes id_{A})-\nu(id_{A} \otimes \nu).\qedhere
\end{displaymath}
\end{proof}
\end{prop}
\begin{rmk}
If $(A,\mu)$ is an associative $\field{K}$-algebra and $id_{A}:A \rightarrow A$ denotes the identity map on $A$, then
\begin{displaymath}
\delta_{H}(id_{A})=\mu(id_{A} \otimes id_{A})+\mu(id_{A} \otimes id_{A})-id_{A}(\mu)=\mu
\end{displaymath}
\end{rmk}
\begin{dfn}[Multiderivations]
The space of the multiderivations on the associative unital $\field{K}$-algebra $(A,\mu,e)$, denoted by $MDer(A)$, is the subalgebra of $(C^{\bullet}(A,A),\smile)$ generated by $Der(A)$.\label{multiderivations}
\end{dfn}

Note that $MDer(A)$ is a graded algebra with
\begin{displaymath}
MDer(A) = \bigoplus_{n \geq 1}MDer^{n}(A),
\end{displaymath}
where $MDer^{n}(A)=MDer(A) \cap C^{n}(A,A)$.
\begin{thm}[The $MDer(A)$ subcomplex]
Every multiderivation is a Hochschild cocycle.
\begin{proof}
Lets proceed by induction to prove that $\delta_{H}$ is identically null on $MDer(A)$. Let $X \in Der(A)$. For all $a,b \in A$
\begin{displaymath}
\delta_{H} X(a \otimes b)=\mu(a \otimes X(b))-X(\mu(a \otimes b))+\mu(X(a) \otimes b)=0
\end{displaymath}
Now, suppose the result for elements in $MDer^{n-1}(A)$ and consider $D \in MDer^{n}(A)$. The space $MDer(A)$ is generated by $Der(A)$, so $D$ can be written as linear combinations of elements of the form $X \smile \tilde{D}$, where $X \in Der(A)$ and $\tilde{D} \in MDer^{n-1}(A)$. By linearity, its enough to consider the evaluation of $\delta_{H}$ in such elements. It follows from the fact that $\delta_{H}$ is a degree 1 derivation on $(C^{\bullet}(A,A),\smile)$ that
\begin{displaymath}
\delta_{H}(X \smile \tilde{D})=(\delta_{H}X)\smile \tilde{D}-X \smile \delta_{H}\tilde{D}=0
\end{displaymath}
Hence, $MDer(A)$ is a subcomplex of $(C^{\bullet}(A,A),\smile)$ and $\delta_{H}$ is identically null on $MDer(A)$.
\end{proof}\label{thm_MDer_trivial_subcomplex}
\end{thm}

The next theorem relates contravariant tensor fields on a differentiable manifold $\manifold{M}$ with multiderivations on the algebra $\funcman{\manifold{M}}$. Before stating the results it is worthy to relate such fields with multiderivations at a point, an analogous to the concept of derivation at a point (see section \ref{algebras}). Lets precise the notion of multiderivation at a point.

Let $\manifold{M}$ be a $m$-dimensional differentiable manifold. The tangent space at every point of $\manifold{M}$, being finite dimensional, allows identify $(T_{p}\manifold{M})^{\otimes n}$ and $(I_{p}/I_{p}^{2})^{\ast \otimes n}$, where $I_{p}$ denotes the ideal of germs of functions which vanish at $p$.

Consider the differentiable tensor bundle $(T\manifold{M})^{\otimes n}$. Let $p \in \manifold{M}$ and $\tau_{p} \in (T\manifold{M})^{\otimes n}$ such that $\tau_{p} \in V_{p}^{\otimes n}=(I_{p}/I_{p}^{2})^{\ast \otimes n}$. By denoting $\mathcal{F}_{p}$ the $\R$-vector space of germs of functions at the point $p$, define the linear map $\vartheta_{p}:\mathcal{F}_{p}^{\otimes n} \rightarrow \R$ given by
\begin{equation*}
\vartheta_{p}(f_{1}\otimes \ldots \otimes f_{n})=
\begin{cases}
	0 & \text{ , if } \exists \; c \in f_{i} \;| \; c(x)=c\;\forall x \in \manifold{M},\\ & \text{ for any }i,\; i=1,\ldots,n\\
	\tau_{p}([f_{1}]\otimes \ldots \otimes [f_{n}]) & \text{ , if } f_{i} \in I_{p},\;\forall i=1,\ldots,n
\end{cases}
\end{equation*}
where $[f_{i}]$ denotes the equivalence class of the germ $f_{i}$ in $I_{p}/I_{p}^{2}$. Since every germ $f$ can be written as $f=\tilde{f}+f(p)$, where $\tilde{f} \in I_{p}$ and $f(p)$ is the germ of the constant function whose value is $f(p)$, $\vartheta_{p}$ satisfies the following:
\begin{eqnarray*}
& & \vartheta_{p}(f_{1}\otimes \ldots \otimes f_{i}g_{i} \otimes \ldots \otimes f_{n}) = \\
& = & g_{i}(p)\vartheta_{p}(f_{1}\otimes \ldots \otimes f_{i} \otimes \ldots \otimes f_{n}) + f_{i}(p)\vartheta_{p}(f_{1}\otimes \ldots \otimes g_{i} \otimes \ldots \otimes f_{n})
\end{eqnarray*}
for every $i$. A linear map $\omega:\mathcal{F}_{p}^{\otimes n} \rightarrow \R$ such that
\begin{eqnarray*}
& & \omega(f_{1}\otimes \ldots \otimes f_{i}g_{i} \otimes \ldots \otimes f_{n}) = \\
& = & g_{i}(p)\omega(f_{1}\otimes \ldots \otimes f_{i} \otimes \ldots \otimes f_{n}) + f_{i}(p)\omega(f_{1}\otimes \ldots \otimes g_{i} \otimes \ldots \otimes f_{n})
\end{eqnarray*}
for every $i=1,\ldots,n$ is called a $\emph{multiderivation of degree n, at the point p}$.

On the other hand, if $\omega:\mathcal{F}_{p}^{\otimes n} \rightarrow \R$ is a multiderivation at $p$, one can relate this to an unique element $\eta_{p} \in V_{p}^{\otimes n}$, such that $\eta_{p}([f_{1}]\otimes \ldots \otimes [f_{n}])=\omega(f_{1}\otimes \ldots \otimes f_{n})$, for all $f_{1}\otimes \ldots \otimes f_{n} \in \mathcal{F}_{p}^{\otimes n}$. To see this, first note that if $c$ represents a constant function, then
\begin{eqnarray*}
& & \omega(f_{1}\otimes \ldots \otimes c \otimes \ldots \otimes f_{n}) = c\; \omega(f_{1}\otimes \ldots \otimes 1 \cdot 1 \otimes \ldots \otimes f_{n}) = \\
& = & c(\omega(f_{1}\otimes \ldots \otimes 1 \otimes \ldots \otimes f_{n})+\omega(f_{1}\otimes \ldots \otimes 1 \otimes \ldots \otimes f_{n})) = \\
& = & 2\;\omega(f_{1}\otimes \ldots \otimes c \otimes \ldots \otimes f_{n})
\end{eqnarray*}
Hence, $\omega(f_{1}\otimes \ldots \otimes c \otimes \ldots \otimes f_{n})=0$. It follows that if $f_{1}\otimes \ldots \otimes f_{n} \in \mathcal{F}_{p}^{\otimes n}$, writing every $f_{i}$ as $f_{i}=\tilde{f}_{i}+f_{i}(p)$, with $\tilde{f}_{i} \in I_{p}$, by linearity of $\omega$, we have
\begin{displaymath}
\omega(f_{1}\otimes \ldots \otimes f_{n})=\omega(\tilde{f}_{1}\otimes \ldots \otimes \tilde{f}_{n})
\end{displaymath}
which shows that the value of $\omega$ is determined only by its value in $I_{p}^{\otimes n}$.

Now, suppose $f_{i} \in I_{p}^{2}$ for some $i$. Then there exists $g_{i},h_{i} \in I_{p}$ such that $f_{i}=g_{i}h_{i}$, resulting
\begin{eqnarray*}
& & \omega(f_{1}\otimes \ldots \otimes f_{i} \otimes \ldots \otimes f_{n}) = \omega(f_{1}\otimes \ldots \otimes g_{i}h_{i} \otimes \ldots \otimes f_{n}) = \\
& = & g_{i}(p)\omega(f_{1}\otimes \ldots \otimes h_{i} \otimes \ldots \otimes f_{n}) + h_{i}(p)\omega(f_{1}\otimes \ldots \otimes g_{i} \otimes \ldots \otimes f_{n}) = 0
\end{eqnarray*}
So, if $\tilde{f}_{i} \equiv \tilde{g}_{i}\;mod\;I_{p}^{2},\;\forall i=1,\ldots,n$, then $\omega(f_{1}\otimes \ldots \otimes f_{n})=\omega(g_{1}\otimes \ldots \otimes g_{n})$ and $\omega$ induces an unique linear map $\eta_{p} \in (I_{p}/I_{p}^{2})^{\ast \otimes n}$. Therefore there exists a bijection between $V_{p}^{\otimes n}$ and multiderivations of degree $n$ at $p$. Hence, an element $\eta_{p} \in V_{p}^{\otimes n}$ can be thought as a multiderivation of degree $n$ at $p$.

Define the action of $\vartheta_{p} \in V_{p}^{\otimes n}$ on elements in $(\funcman{\manifold{M}})^{\otimes n}$ by $\vartheta_{p}(F_{1}\otimes \ldots \otimes F_{n})=\vartheta_{p}(f_{1}\otimes \ldots \otimes f_{n})$ on decomposable elements and extended it by linearity, where $f_{i}$ denotes a representing element of $F_{i} \in \funcman{\manifold{M}}$ in $\mathcal{F}_{p}$, $i=1,\ldots,n$.

We can now prove the following
\begin{thm}
Let $\manifold{M}$ be an $m$-dimensional differentiable manifold. Then
\begin{displaymath}
\Gamma((T\manifold{M})^{\otimes n}) \approx_{\catvect{R}} MDer^{n}(\funcman{\manifold{M}}),\quad \forall n \geq 1.
\end{displaymath}
\begin{proof}
Given $\tau \in \Gamma((T\manifold{M})^{\otimes n})$, define a linear map $\bar{\tau}:\funcman{\manifold{M}}^{\otimes n}\rightarrow \funcman{\manifold{M}}$, given by
\begin{displaymath}
\bar{\tau}(f_{1}\otimes \ldots \otimes f_{n})(p)=\tau_{p}(f_{1}\otimes \ldots \otimes f_{n}),\;\forall p \in \manifold{M}
\end{displaymath}
Note that $\tau_{p} \in (I_{p}/I_{p}^{2})^{\ast \otimes n}$ means that $\tau_{p}$ can be written as linear combination of elements of the form $v^{1}_{p}\otimes \ldots \otimes v^{n}_{p}$ with $v^{i}_{p} \in (I_{p}/I_{p}^{2})^{\ast},\;\forall i=1,\ldots,n$, and $(v^{1}_{p}\otimes \ldots \otimes v^{n}_{p})(f_{1}\otimes \ldots \otimes f_{n})=v^{1}_{p}(f_{1}) \ldots v^{n}_{p}(f_{n})$. Since $\tau$ is a differentiable section, we have $\bar{\tau} \in MDer^{n}(\funcman{\manifold{M}})$.

The assignment $\tau \mapsto \bar{\tau}$ is clearly linear. We claim it is a bijection.

To show that it is injective, it is enough to see that
\begin{eqnarray*}
\bar{\tau}(f_{1}\otimes \ldots \otimes f_{n}) = 0 & , & \forall f_{i} \in \funcman{\manifold{M}},\; i=1,\ldots,n \Rightarrow \\
\Rightarrow \bar{\tau}(f_{1}\otimes \ldots \otimes f_{n})(p) = 0 & , & \forall p \in \manifold{M},\;\forall f_{i} \in \funcman{\manifold{M}},\; i=1,\ldots,n \Rightarrow \\
\Rightarrow \tau_{p}(f_{1}\otimes \ldots \otimes f_{n}) = 0 & , & \forall p \in \manifold{M},\;\forall f_{i} \in \funcman{\manifold{M}},\; i=1,\ldots,n \Rightarrow \\
\Rightarrow \tau_{p}=0 & , & \forall p \in \manifold{M} \Rightarrow \\
\Rightarrow \tau=0
\end{eqnarray*}
where in the last but one step it was used the following fact. Taking a local chart $(U,\varphi)$ at the point $p$, with $\varphi(p)=0$ and taking $x^{i}=t^{i}\circ \varphi$, with $t^{i}:\R^{m}\rightarrow \R$ the canonical projection on the $i$-th component (see section \ref{algebras}), we can write
\begin{displaymath}
\tau_{p}(f_{1}\otimes \ldots \otimes f_{n})=\sum_{(i_{1},\ldots,i_{n})}a^{i_{1}\ldots i_{n}}\frac{\partial}{\partial x^{i_{1}}}\otimes \ldots \otimes \frac{\partial}{\partial x^{i_{n}}}
\end{displaymath}
where $(i_{1},\ldots,i_{n})$ under the summation symbol means that the sum must be evaluated for each $i_{j}$, with $j=1,\ldots,n$, $i_{j}=1,\ldots,m$. Evaluating $\tau_{p}$ on elements of the form $(x^{i_{1}}\otimes \ldots \otimes x^{i_{n}})$ successively, we have $a^{i_{1}\ldots i_{n}}=0$ for any combination of indices $i_{j}$.

To show that it is surjective, consider $D \in MDer^{n}(\funcman{\manifold{M}})$. For each $p \in \manifold{M}$, define the linear map $\tau_{p}:\mathcal{F}_{p}^{\otimes n}\rightarrow \R$ given by
\begin{displaymath}
\tau_{p}(f_{1}\otimes \ldots \otimes f_{n})=D(f_{1}\otimes \ldots \otimes f_{n})(p),\;\forall f_{1}\otimes \ldots \otimes f_{n} \in \mathcal{F}_{p}
\end{displaymath}
where $f_{i}$ on the left hand side denotes the germ of the function $f_{i}$ written on the right. For now on in this proof we will denote germs and functions by the same symbol to avoid cumbersome notation. It is clear when a symbol is for a function or for a germ from the operator on such symbol. Lets show that $\tau_{p}$ is well defined. Let $f_{i},g_{i} \in \funcman{\manifold{M}}$ with $i=1,\ldots,n$, such that $f_{i} \equiv g_{i}$ on $\mathcal{F}_{p}$, for each $i=1,\ldots,n$. Let $W_{i}$ be open neighbourhoods of $p \in \manifold{M}$ such that $f_{i}|_{W_{i}}=g_{i}|_{W_{i}},\;i=1,\ldots,n$. Let $(V,\varphi)$ be a local chart such that $\varphi(p)=0$. Taking $U=V \cap W_{1} \cap \ldots \cap W_{n}$ we have $(U,\varphi)$ still a local chart around $p$. If necessary, we can shrink $U$ to make $\varphi(U)$ open and convex on $\R^{m}$. As $f_{i}$ and $g_{i}$ coincides on $U$ for each $i$, it also coincides $\tilde{f}_{i}=f_{i}-f_{i}(p)$ and $\tilde{g}_{i}=g_{i}-g_{i}(p)$ on $U$ for each $i$. Hence,
\begin{displaymath}
\frac{\partial (\tilde{f}_{i} \circ \varphi^{-1})}{\partial t^{j}}\biggl |_{0}=\frac{\partial (\tilde{g}_{i} \circ \varphi^{-1})}{\partial t^{j}}\biggl |_{0},\;\forall j=1,\ldots,m,\;\forall i=1,\ldots,n.
\end{displaymath}

By $D \in MDer^{n}(\funcman{\manifold{M}})$, we have
\begin{displaymath}
D(f_{1} \otimes \ldots \otimes f_{n})=D((\tilde{f}_{1}+f_{1}(p)) \otimes \ldots \otimes (\tilde{f}_{n}+f_{n}(p)))=D(\tilde{f}_{1} \otimes \ldots \otimes \tilde{f}_{n})
\end{displaymath}
which results in
\begin{eqnarray*}
& & \tau_{p}(f_{1} \otimes \ldots \otimes f_{n}) = D(f_{1} \otimes \ldots \otimes f_{n})(p) = D(\tilde{f}_{1} \otimes \ldots \otimes \tilde{f}_{n})(p) = \\
& & = \sum_{j_{1},\ldots,j_{n}=1}^{m}\frac{\partial (\tilde{f}_{i} \circ \varphi^{-1})}{\partial t^{j_{1}}}\biggl |_{0} \ldots \frac{\partial (\tilde{f}_{i} \circ \varphi^{-1})}{\partial t^{j_{n}}}\biggl |_{0} D(x^{j_{1}} \otimes \ldots \otimes x^{j_{n}})(p) + \\
& & + \sum_{\substack{j_{1},\ldots,j_{n} = 1 \\ l_{1},\ldots,l_{n} = 1}}^{m} \int_{0}^{1}(1-s)\frac{\partial^{2}(\tilde{f}_{1} \circ \varphi^{-1})}{\partial t^{j_{1}} \partial t^{l_{1}}}\biggl |_{sy}ds \ldots \int_{0}^{1}(1-s)\frac{\partial^{2}(\tilde{f}_{n} \circ \varphi^{-1})}{\partial t^{j_{n}} \partial t^{l_{n}}}\biggl |_{sy}ds \cdot \\
& & \cdot D(x^{j_{1}}x^{l_{1}} \otimes \ldots \otimes x^{j_{n}}x^{l_{n}})(p) = \\
& & = \sum_{j_{1},\ldots,j_{n}=1}^{m}\frac{\partial (\tilde{g}_{i} \circ \varphi^{-1})}{\partial t^{j_{1}}}\biggl |_{0} \ldots \frac{\partial (\tilde{g}_{i} \circ \varphi^{-1})}{\partial t^{j_{n}}}\biggl |_{0} D(x^{j_{1}} \otimes \ldots \otimes x^{j_{n}})(p) + \\
& & + \sum_{\substack{j_{1},\ldots,j_{n} = 1 \\ l_{1},\ldots,l_{n} = 1}}^{m} \int_{0}^{1}(1-s)\frac{\partial^{2}(\tilde{g}_{1} \circ \varphi^{-1})}{\partial t^{j_{1}} \partial t^{l_{1}}}\biggl |_{sy}ds \ldots \int_{0}^{1}(1-s)\frac{\partial^{2}(\tilde{g}_{n} \circ \varphi^{-1})}{\partial t^{j_{n}} \partial t^{l_{n}}}\biggl |_{sy}ds \cdot \\
& & \cdot D(x^{j_{1}}x^{l_{1}} \otimes \ldots \otimes x^{j_{n}}x^{l_{n}})(p) = D(\tilde{g}_{1} \otimes \ldots \otimes \tilde{g}_{n})(p) = \\
& & = D(g_{1} \otimes \ldots \otimes g_{n})(p) = \tau_{p}(g_{1} \otimes \ldots \otimes g_{n})
\end{eqnarray*}
because $D(x^{j_{1}}x^{l_{1}} \otimes \ldots \otimes x^{j_{n}}x^{l_{n}})(p)=0$ for any combination of indices $(j_{k},l_{k})$. Thus, $\tau_{p}$ is well defined as a linear map on $\mathcal{F}_{p}^{\otimes n}$ to real values, for all $p \in \manifold{M}$. Besides that, we have
\begin{multline*}
\tau_{p}(f_{1} \otimes \ldots \otimes f_{i}g_{i} \otimes \ldots \otimes f_{n}) = D(f_{1} \otimes \ldots \otimes f_{i}g_{i} \otimes \ldots \otimes f_{n})(p) = \\
= g_{i}(p)D(f_{1} \otimes \ldots \otimes f_{i} \otimes \ldots \otimes f_{n})(p) + f_{i}(p)D(f_{1} \otimes \ldots \otimes g_{i} \otimes \ldots \otimes f_{n})(p) = \\
= g_{i}(p)\tau_{p}(f_{1} \otimes \ldots \otimes f_{i} \otimes \ldots \otimes f_{n}) + f_{i}(p)\tau_{p}(f_{1} \otimes \ldots \otimes g_{i} \otimes \ldots \otimes f_{n})
\end{multline*}
on each entry. Therefore $\tau_{p}$ is a multiderivation of degree $n$ at the point $p$, for each $p \in \manifold{M}$. Finally, lets construct the map $\tau:\manifold{M}\rightarrow (T\manifold{M})^{\otimes n}$ given by $\tau(p)=\tau_{p}$. The map $\tau$ is a differentiable section, because for each $p \in \manifold{M}$
\begin{displaymath}
\tau(p)(f_{1} \otimes \ldots \otimes f_{n}) = \tau_{p}(f_{1} \otimes \ldots \otimes f_{n}) = D(f_{1} \otimes \ldots \otimes f_{n})(p)
\end{displaymath}
and $D(f_{1} \otimes \ldots \otimes f_{n}) \in \funcman{\manifold{M}}$ for any linear combination of elements $f_{1} \otimes \ldots \otimes f_{n} \in \funcman{\manifold{M}}$.

Thus, the stated assignment is an isomorphism of $\R$-vector spaces between $\Gamma((T\manifold{M})^{\otimes n})$ and $MDer^{n}(\funcman{\manifold{M}})$.
\end{proof}\label{isomorphism_tensor_fields_MDer}
\end{thm}
The last theorem reveals that if $\manifold{M}$ is a $m$-dimensional differentiable manifold, an element in $D \in MDer^{n}(\funcman{\manifold{M}})$ can be written in local coordinates as
\begin{displaymath}
D = \sum_{j_{1},\ldots,j_{n}=1}^{m}D(x^{j_{1}} \otimes \ldots \otimes x^{j_{n}})\frac{\partial}{\partial x^{j_{1}}}\otimes \ldots \otimes \frac{\partial}{\partial x^{j_{n}}}.
\end{displaymath}
\section{Iterated Derivations}
\begin{dfn}[Iterated Derivation]
Let $A$ be a commutative associative unital $\field{K}$-algebra. The space of iterated derivations, denoted by $SDer(A)$, is the subalgebra of $(C^{1}(A,A),\circ)$ generated by $Der(A)$. We denote by $SDer^{n}(A)$ the set of elements $D \in SDer(A)$ which can be written as linear combinations of elements of the form $X_{1} \circ \ldots \circ X_{r}$, with $X_{i} \in Der(A), \;\forall i=1,\ldots,r,\;r \leq n$.\label{SDer}
\end{dfn}
\begin{rmk}
Note that $SDer(A)$ can not be written as direct sum of the spaces $SDer^{n}(A)$. However, if $r \leq n$ then we have $SDer^{r}(A) \subset SDer^{n}(A)$. Hence, we have a filtration on the algebra $SDer(A)$.
\end{rmk}
\begin{thm}
If $D \in SDer^{n}(A)$, then $D$ is a derivation of order $\leq$ $n$.\footnote{A similar notion for Lie algebras can be found in \cite{sardanashvily_diff_op_Lie_algebras}.}
\begin{proof}
Denote the product on $A$ by juxtaposition. We proceed by induction on $n$. Surely, if $X \in Der(A)$, then $X$ is a derivation of order $\leq$ 1. Suppose $D \in SDer^{n}(A)$ and the result valid for $n-1$. By linearity, it is enough to consider $D$ as $D=\tilde{D} \circ X_{n}$, where $\tilde{D} \in SDer^{n-1}(A)$ and $X_{n} \in Der(A)$. By the induction hypothesis and the fact that $SDer^{r-1}(A) \subset SDer^{r}(A)$ for all $r \geq 1$, it is enough to consider the high order terms of $\tilde{D}$, \emph{i.e.} terms as $X_{1} \circ \ldots \circ X_{n-1}$. To show that $X_{1} \circ \ldots \circ X_{n-1}$ is a differential operator of order $\leq$ $n$, given $a \in A$, we must show that the operator $\Delta_{a}$, given by
\begin{displaymath}
\Delta_{a}(b)=(X_{1} \circ \ldots \circ X_{n})(ab)-a(X_{1} \circ \ldots \circ X_{n})(b)
\end{displaymath}
for all $b \in A$, is a differential operator of order $\leq$ $n-1$. We have
\begin{eqnarray*}
& & \Delta_{a}(b) = (X_{1} \circ \ldots \circ X_{n})(ab)-a(X_{1} \circ \ldots \circ X_{n})(b) = \\
& & = (X_{1} \circ \ldots \circ X_{n})(a) \cdot b + \sum_{i=1}^{n}(X_{1} \circ \ldots \circ \hat{X}_{i} \circ \ldots \circ X_{n})(a)X_{i}(b) + \\
& & + \sum_{1 \leq i < j \leq n}(X_{1} \circ \ldots \circ \hat{X}_{i} \circ \ldots \circ \hat{X}_{j} \circ \ldots \circ  X_{n})(a)(X_{i} \circ X_{j})(b) + \ldots + \\
& & + \sum_{I_{k}}X_{\hat{I}_{k}}(a)X_{I_{k}}(b) + \ldots + \sum_{i=1}^{n}X_{i}(a)(X_{1} \circ \ldots \circ \hat{X}_{i} \circ \ldots \circ X_{n})(b)
\end{eqnarray*}
where $I_{k}$ represents a set of indices, subset of $\{ 1,\ldots,n \}$, with exactly $k$ elements $\{ i_{1},\ldots,i_{k} \}$, such that $i_{1} < \ldots < i_{k}$, $X_{\hat{I}_{k}}$ represents the composition $X_{1} \circ \ldots \circ \hat{X}_{i_{j}} \circ \ldots \circ X_{n}$ in which are absent all elements $X_{i_{1}},\ldots,X_{i_{k}}$ in this order, and $X_{I_{k}}$ represents the composition $X_{i_{1}} \circ \ldots \circ X_{i_{k}}$. Hence, we have $\Delta_{a}$ operator acting on $b$ with composites having at most $n-1$ factors. Therefore $\Delta_{a} \in SDer^{n-1}(A)$, which is by the induction hypothesis a differential operator of order $\leq$ $n-1$, for all $a \in A$. Thus $D$ is a differential operator of order $\leq$ $n$. By considering operators as $\tilde{D} \circ X$, with $X \in Der(A)$, it is clear that $\tilde{D}(X(\alpha))=0$, for all $\alpha \in \field{K}$ (properly identified as element of $A$). Hence, $D$ is a derivation of order $\leq$ $n$.
\end{proof}
\end{thm}
\begin{thm}
Let $\manifold{M}$ be an $m$-dimensional differentiable manifold. If $D$ is a derivation of order $\leq$ $r$ on $\funcman{\manifold{M}}$, then $D \in SDer^{r}(\funcman{\manifold{M}})$.

\begin{proof}
We proceed by induction. If $D$ is a derivation of order $\leq$ 1, then $D \in Der(\funcman{\manifold{M}})$ therefore $D \in SDer^{1}(\funcman{\manifold{M}})$. Suppose the result for $r-1$. Let $D$ be a derivation of order $\leq$ $r$ on $\funcman{\manifold{M}}$. Then, by theorem \ref{iso_deriv_order_r_jets}, $D$ can be related to an element $D \in \Gamma(J^{r}(\manifold{M}))$. For each $p \in \manifold{M}$, define the linear map $\Phi_{r,p}:I_{p}/I_{p}^{r+1} \rightarrow I_{p}/I_{p}^{r}$ which associates the equivalence class of a germ of a function $f$ in $I_{p}/I_{p}^{r+1}$ to its class in $I_{p}/I_{p}^{r}$. This is well defined, because $I_{p}^{r+1} \subset I_{p}^{r}$ and it is a projection because, by the Taylor's formula, if $f$ has a representing in $I_{p}/I_{p}^{r}$, then it has a representing in $I_{p}/I_{p}^{r+1}$ such that $[f]_{r}=\Phi_{r,p}([f]_{r+1})$. Note that if $f \in I_{p}^{r} \;mod\; I_{p}^{r+1}$, then $\Phi_{r,p}(f) = 0$, and by the other hand, if $\Phi_{r,p}(f) = 0$, then $f \in I_{p}^{r}\;mod\;I_{p}^{r+1}$. Thus, $Ker(\Phi_{r,p}) \approx_{\catvect{R}} I_{p}^{r}/I_{p}^{r+1}$. We have, naturally, $I_{p}/I_{p}^{r+1} \approx_{\catvect{R}} I_{p}/I_{p}^{r} \oplus I_{p}^{r}/I_{p}^{r+1}$.

The dual map to $\Phi_{r,p}$ is $\Phi_{r,p}^{\ast}:J^{r-1}_{p} \rightarrow J^{r}_{p}$, given by
\begin{displaymath}
(\Phi_{r,p}^{\ast}(u))(f)=u(\Phi_{r,p}(f))
\end{displaymath}
remembering that $J^{r}_{p}=(I_{p}/I_{p}^{r+1})^{\ast}$. $\Phi_{r,p}^{\ast}$ is injective. This follows from the fact of being dual to a surjective linear map between vector spaces, because if $u \in J^{r-1}_{p}$ is such that $\Phi_{r,p}^{\ast}(u) = 0$, then $(\Phi_{r,p}^{\ast}(u))(f) = 0$ for all $f \in I_{p}/I_{p}^{r+1}$ and then, $u(\Phi_{r,p}(f)) = 0$ for all $f \in I_{p}/I_{p}^{r+1}$. As $\Phi_{r,p}$ is surjective, given $g \in I_{p}/I_{p}^{r}$, there exists $f \in I_{p}/I_{p}^{r+1}$ such that $g = \Phi_{r,p}(f)$. Hence, $u(g) = 0$ for all $g \in I_{p}/I_{p}^{r}$ therefore $u=0$.

The map $\Phi_{r}^{\ast}:J^{r-1}(\manifold{M}) \rightarrow J^{r}(\manifold{M})$ such that $\Phi_{r}^{\ast}(\xi) = \Phi_{r,\pi(\xi)}^{\ast}(\xi)$ is a morphism of differentiable vector bundles. We have $\Phi_{r}^{\ast}$ fibre preserving and linear on fibres by construction. Furthermore, if $\xi \in J^{r-1}(\manifold{M})$, locally, $\xi$ is written as
\begin{displaymath}
\xi = \sum_{k=1}^{r-1}\sum_{1 \leq i_{1} \leq \ldots \leq i_{k} \leq m}\xi(x^{i_{1}} \ldots x^{i_{k}})\frac{\partial^{k}}{\partial x^{i_{1}} \ldots \partial x^{i_{k}}}
\end{displaymath}
But $\Phi_{r}^{\ast}(\xi)$ is written locally as
\begin{displaymath}
\xi = \sum_{k=1}^{r-1}\sum_{1 \leq i_{1} \leq \ldots \leq i_{k} \leq m}\xi(y^{i_{1}} \ldots y^{i_{k}})\frac{\partial^{k}}{\partial y^{i_{1}} \ldots \partial y^{i_{k}}}
\end{displaymath}
because terms of order $r$ does not belong to the range of $\Phi_{r}^{\ast}$. By the fact that local charts on $J^{r-1}(\manifold{M})$ and $J^{r}(\manifold{M})$ are fibred charts, there exists a diffeomorphism sending the coordinate expression of $\xi$ in terms of $y^{i}$ and derivatives, to the coordinate expression of $\xi$ in terms of $x^{i}$ and derivatives. By the match of those expressions follows $\Phi_{r}^{\ast}$ differentiable.

Given $p \in \manifold{M}$, by the induction hypothesis and the inclusion above, it is enough to consider derivations of order $\leq$ $r$ such that in a neighbourhood of $p$ have only terms of order $r$. Let $\eta$ a such derivation and $(U,x^{1},\ldots,x^{m})$ a local chart around $p$ for which this occurs. $\eta$ being a derivation of order $\leq$ $r$ leads to
\begin{equation*}
\eta(x^{i_{1}}\ldots x^{i_{r}}) = \Delta_{i_{r}}(x^{i_{1}} \ldots x^{i_{r-1}}) + x^{i_{r}}\eta(x^{i_{1}}\ldots x^{i_{r-1}})
\end{equation*}
with $\Delta_{i_{r}}$ differential operator of order $\leq$ $r-1$. By the choice of the local chart, we have $\eta(x^{i_{1}}\ldots x^{i_{r-1}}) = 0$, because $\eta$ has only terms of order $r$. Therefore
\begin{displaymath}
\eta(x^{i_{1}}\ldots x^{i_{r}}) = \Delta_{i_{r}}(x^{i_{1}} \ldots x^{i_{r-1}})
\end{displaymath}
and from this follows that $\Delta_{i_{r}}$ is a derivation of order $\leq$ $r-1$. $\eta$ can be written in terms of $\Delta_{i_{r}}$ in this way
\begin{eqnarray}
\eta & = & \sum_{k=1}^{m}\sum_{i_{1} \leq \ldots \leq i_{r-1}}\frac{\eta(x^{i_{1}}\ldots x^{i_{r-1}}x^{k})}{r!}\frac{\partial^{r}}{\partial x^{i_{1}} \ldots \partial x^{i_{r-1}} \partial x^{k}} = \nonumber \\
& = & \sum_{k=1}^{m}\sum_{i_{1} \leq \ldots \leq i_{r-1}}\frac{\Delta_{k}(x^{i_{1}}\ldots x^{i_{r-1}})}{r!}\frac{\partial^{r-1}}{\partial x^{i_{1}} \ldots \partial x^{i_{r-1}}}\frac{\partial}{\partial x^{k}} = \nonumber \\
& = & \sum_{k=1}^{m}\frac{\Delta_{k}}{r!}\frac{\partial}{\partial x^{k}}
\label{eq_op_dif_escrito_composta}
\end{eqnarray}
As $\Delta_{k}$ is a derivation of order $\leq$ $r-1$, the induction hypothesis allows to write
\begin{displaymath}
\Delta_{k} = v_{k} \circ u_{k}
\end{displaymath}
where $v_{k}$ is a vector field and $u_{k}$ is a derivation of order $\leq$ $r-2$, both defined on $U$.

By the equation \ref{eq_op_dif_escrito_composta}, we have
\begin{equation*}
\eta = \sum_{k=1}^{m}\frac{\Delta_{k}}{r!}\frac{\partial}{\partial x^{k}} = \sum_{k=1}^{m}\frac{(v_{k} \circ u_{k})}{r!}\frac{\partial}{\partial x^{k}} = \sum_{k=1}^{m}v_{k} \left ( \frac{u_{k}}{r!}\frac{\partial}{\partial x^{k}} \right )
\end{equation*}
as the term $\frac{u_{k}}{r!}\frac{\partial}{\partial x^{k}}$ is a composition of derivations, it is itself a derivation of order $\leq$ $r-1$ defined on $U$. Therefore
\begin{displaymath}
\eta = \sum_{k=1}^{m}v_{k} \circ w_{k}
\end{displaymath}
with $v_{k} \in \Gamma(J^{1}(U))$ and $w_{k} \in \Gamma(J^{r-1}(U))$, for each $k=1,\ldots,m$. For the sake of simplicity, we denote this by $\eta = v \circ u$.

Let $\{ U_{\alpha} \}$ be a locally finite open covering of $\manifold{M}$ and $\{ \rho_{\alpha} \}$ a partition of unity subordinated to such covering. For each index $\alpha$, we can find $v_{\alpha}$ and $u_{\alpha}$ as above, such that
\begin{displaymath}
\eta = v_{\alpha} \circ u_{\alpha}
\end{displaymath}
Lets construct the fields $\zeta \in \mathfrak{X}(\manifold{M})$, $\xi,\theta \in \Gamma(J^{r-1}(\manifold{M}))$ by
\begin{equation*}
\zeta = \sum_{\lambda}\rho_{\lambda}v_{\lambda}, \quad \xi = \sum_{\nu}\rho_{\nu}u_{\nu}, \quad \theta = \sum_{\beta}\gamma_{\beta}u_{\beta}
\end{equation*}
where $\gamma_{\beta} = \sum_{\alpha}\rho_{\alpha}v_{\alpha}(\rho_{\beta})$. Note that $\theta$ is well defined, because if $\rho_{\beta}$ has support on $U_{\beta}$, so are its derivatives and then, given $p \in \manifold{M}$, $\gamma_{\beta}(p)$ does not vanish only for a finite number of indices $\beta$. Furthermore, given $f \in \funcman{\manifold{M}}$
\begin{displaymath}
\rho_{\lambda}v_{\lambda}(\rho_{\nu}u_{\nu}(f)) = \rho_{\lambda}\rho_{\nu}v_{\lambda}(u_{\nu}(f)) + \rho_{\lambda}v_{\lambda}(\rho_{\nu})u_{\nu}(f)
\end{displaymath}
leading to
\begin{displaymath}
\rho_{\lambda}\rho_{\nu}\eta(f) = \rho_{\lambda}v_{\lambda}(\rho_{\nu}u_{\nu}(f)) - \rho_{\lambda}v_{\lambda}(\rho_{\nu})u_{\nu}(f)
\end{displaymath}
because if $U_{\lambda} \cap U_{\nu} = \varnothing$, then either $\rho_{\lambda}$ or $\rho_{\nu}$ vanish, and then $\rho_{\lambda}\rho_{\nu}v_{\lambda} \circ u_{\nu} = \rho_{\lambda}\rho_{\nu}\eta$, and if $U_{\lambda} \cap U_{\nu} \neq \varnothing$, then $u_{\nu p} = u_{\lambda p}$ at each $p \in U_{\lambda} \cap U_{\nu}$, leading to $\rho_{\lambda}\rho_{\nu}v_{\lambda} \circ u_{\nu} = \rho_{\lambda}\rho_{\nu}\eta$.

Hence,
\begin{eqnarray*}
\eta(f) & = & \sum_{\lambda,\nu}\rho_{\lambda}\rho_{\nu}\eta(f) = \sum_{\lambda,\nu}\rho_{\lambda}v_{\lambda}(\rho_{\nu}u_{\nu}(f)) - \sum_{\lambda,\nu}\rho_{\lambda}v_{\lambda}(\rho_{\nu})u_{\nu}(f) = \\
& = & \sum_{\lambda}\rho_{\lambda}v_{\lambda} \left ( \sum_{\nu}\rho_{\nu}u_{\nu}(f) \right ) - \sum_{\nu}\gamma_{\nu}u_{\nu}(f) = \\
& = & (\zeta \circ \xi)(f) - \theta(f)
\end{eqnarray*}

By the induction hypothesis, $\xi,\theta \in \Gamma(J^{r-1}(\manifold{M})$ can be related to elements in $SDer^{r-1}(\funcman{\manifold{M}})$, leading to $\eta \in SDer^{r}(\funcman{\manifold{M}})$. As an arbitrary element $D \in \Gamma(J^{r}(\manifold{M})$ is a linear combination of elements in $\Gamma(J^{r-1}(\manifold{M})$ and elements of order $r$, it follows that
\begin{displaymath}
D \in SDer^{r}(\funcman{\manifold{M}}). \qedhere
\end{displaymath}
\end{proof}

\label{thm_order_r_deriv_are_elements_in_SDer}
\end{thm}
\section{Polydifferential Operators}
The Hochschild-Kostant-Rosenberg theorem is usually stated as an isomorphism of graded algebras between Hochschild homology and universal differential forms (given by the K\"ahler differentials) of a smooth algebra. A proof of this version can be found in \cite{loday_cyclic_homology}. However, we want a dual version of this fact, by relating Hochschild cohomology of an algebra and its multilinear transformations. The process of taking duals often involves some restriction to a nice subspace. For infinite dimensional cases, the dual of a vector space is too big and an analogous copy of the original space that retains or preserves the desired properties lies in a specific kind of subspace. In the case of the Hochschild-Kostant-Rosenberg theorem we must restrict the Hochschild cohomology to the subcomplex of polyderivations.
\begin{dfn}[Polyderivations on an algebra]
Let $A$ be a commutative associative unital $\field{K}$-algebra. The space of polyderivations on the algebra $A$, denoted by $D_{poly}(A)$, is the subalgebra of $(C^{\bullet}(A,A),\smile)$ generated by $SDer(A)$. We denote $D_{poly}^{n}(A)=D_{poly}(A)\cap C^{n}(A,A)$. Also, we denote by $D_{poly}^{n,r}(A)$ the space of polyderivations of degree $n$ and order $\leq$ $r$ i.e. elements in $C^{n}(A,A)$ which are polyderivations generated by $SDer^{r}(A)$.\label{polyderivations}
\end{dfn}
\begin{thm}
$(D_{poly}(A),\delta_{H})$ is a filtered subcomplex of $(C^{\bullet}(A,A),\delta_{H})$.
\begin{proof}
For the sake of simplicity, we denote the product on $A$ by juxtaposition. Take an element $D \in D_{poly}^{n,r}(A)$. Then $D$ is a linear combination of elements of the form $D_{1} \smile \ldots \smile D_{n}$, with $D_{i} \in SDer^{r}(A)$, for all $i=1,\ldots,n$. However, if $D_{i} \in SDer^{r}(A)$, then it is linear combination of elements of the form $X^{i}_{1} \circ \ldots \circ X^{i}_{j}$, $j \leq r$, with $X^{i}_{j} \in Der(A)$, for all $i=1,\ldots,n$, for all $j \leq r$. Then, if $a,b \in A$
\begin{eqnarray}
& & \delta_{H}(X^{i}_{1} \circ \ldots \circ X^{i}_{j})(a \otimes b) = \nonumber \\
& & = a(X^{i}_{1} \circ \ldots \circ X^{i}_{j})(b) - (X^{i}_{1} \circ \ldots \circ X^{i}_{j})(ab) + (X^{i}_{1} \circ \ldots \circ X^{i}_{j})(a)b = \nonumber \\
& & = -\sum_{k=1}^{j-1}\sum_{I_{k}}(X^{i}_{\hat{I}_{k}})(a) (X^{i}_{I_{k}})(b)
\label{formula_dif_hochschild_on_SDer_r}
\end{eqnarray}
where $I_{k}$ denotes a set of indices, subset of $\{ 1,\ldots,j \}$, with exactly $k$ elements $l_{1},\ldots,l_{k}$ such that $l_{1} < \ldots < l_{k}$, for $k \leq j$, $X^{i}_{\hat{I}_{k}}$ denotes the composite $X^{i}_{1} \circ \ldots \circ \hat{X}^{i}_{l_{s}} \circ \ldots \circ X^{i}_{j}$, in which are absent all elements $X^{i}_{l_{s}}$, $l_{s} \in I_{k}$, in order, and $X^{i}_{I_{k}}$ denotes the composite $X^{i}_{l_{1}} \circ \ldots \circ X^{i}_{l_{k}}$ in that order.

Hence, $\delta_{H}(X^{i}_{1} \circ \ldots \circ X^{i}_{j}) \in D_{poly}^{2,j-1}(A)$. As $\delta_{H}$ is a degree 1 derivation on $(C^{\bullet}(A,A),\smile)$, it follows that
\begin{equation}
\delta_{H}(D_{1} \smile \ldots \smile D_{n}) = \sum_{i=1}^{n}(-1)^{i+1}D_{1} \smile \ldots \smile \delta_{H}(D_{i}) \smile \ldots \smile D_{n}
\label{formula_dif_hochschild_on_D_poly}
\end{equation}
By linearity, $D \in D_{poly}^{n,r}(A)$, results $\delta_{H}(D) \in D_{poly}^{n+1,r}(A)$. It shows that $(D_{poly}(A),\delta_{H})$ is subcomplex of $(C^{\bullet}(A,A),\delta_{H})$, filtered by order of derivations.
\end{proof}
\end{thm}
\begin{dfn}[Alternator on $D_{poly}^{n,r}(A)$]
If $A$ is a commutative associative unital $\field{K}$-algebra, where $\field{K}$ is a field with characteristic 0, we define for $n \geq 1$ the linear map $Alt:D_{poly}^{n,r}(A) \rightarrow D_{poly}^{n,r}(A)$ given, on decomposable elements, by
\begin{displaymath}
Alt(D_{1} \smile \ldots \smile D_{n})=\frac{1}{n!}\sum_{\sigma \in S_{n}}\varepsilon(\sigma)D_{\sigma(1)} \smile \ldots \smile D_{\sigma(n)}
\end{displaymath}
where $\sigma$ denotes a permutation in $S_{n}$, the set of all permutations on $n$ elements, and $\varepsilon(\sigma)$ denotes the signal of this permutation.
\end{dfn}
\begin{prop}\label{cocycle_coboundary_n_vector_hochschild_formula}
Let $D \in D_{poly}^{n,r}(\funcman{\manifold{M}})$ such that $D$ is closed for the Hochschild differential. Then there exists a cochain $E \in D_{poly}^{n-1,r+1}(\funcman{\manifold{M}})$ and an alternating element $\eta \in MDer^{n}(\funcman{\manifold{M}})$ such that
\begin{equation}
D = \delta_{H}(E)+\eta
\end{equation}
The proof of the proposition is quite technical and can be found in \cite{s_gutt_equivalence_star_products}.
\end{prop}
\begin{rmk}
Let $A$ be a commutative associative unital $\field{K}$-algebra. We denote $\mathcal{D}(A) = A \oplus D_{poly}(A)$. Note that $(\mathcal{D}(A),\delta_{H})$ is subcomplex of the Hochschild complex $(C^{\bullet}(A,A),\delta_{H})$.
\end{rmk}
\begin{thm}[The Hochschild-Kostant-Rosenberg theorem for differentiable manifolds\footnote{This proof follows the technique in \cite{s_gutt_local_cohomology}}]
Let $\manifold{M}$ be a $m$-dimensional differentiable manifold. There is a quasi-isomorphism between the complexes $(\mathcal{D}(\funcman{\manifold{M}}),\delta_{H})$ and $(\Omega_{\bullet}(\manifold{M}),d)$, where $d:\Omega_{\bullet}(\manifold{M}) \rightarrow \Omega_{\bullet}(\manifold{M})$ is the null differential on polyvector fields $\Omega_{\bullet}(\manifold{M})=\Gamma(\Lambda T\manifold{M})$.
\begin{proof}
Let $Alt(MDer^{n}(\funcman{\manifold{M}}))$ be the range of the alternator on\linebreak$MDer^{n}(\funcman{\manifold{M}}) = D_{poly}^{n,1}(\funcman{\manifold{M}})$. Define the linear map $\psi:\Omega_{n}(\manifold{M}) \rightarrow Alt(MDer^{n}(\funcman{\manifold{M}}))$ given, on decomposable elements, by
\begin{displaymath}
\psi(X_{1} \wedge \ldots \wedge X_{n}) = Alt(X_{1} \smile \ldots \smile X_{n})
\end{displaymath}
for $n \geq 1$. Note that $\psi$ is fibre preserving. Lets show that $\psi$ is injective. Let $\eta \in \Omega_{n}(\manifold{M})$ such that $\psi(\eta) = 0$. At each $p \in \manifold{M}$, $\eta_{p}$ is written as linear combination of elements in a base for $\Lambda_{p}(T_{p}\manifold{M})$, of the form $X_{i_{1} p} \wedge \ldots \wedge X_{i_{n} p}$. However,
\begin{eqnarray*}
& & \psi(X_{i_{1} p} \wedge \ldots \wedge X_{i_{n} p}) = Alt(X_{i_{1} p} \smile \ldots \smile X_{i_{n} p}) = Alt(X_{i_{1} p} \otimes \ldots \otimes X_{i_{n} p}) = \\
& & = X_{i_{1} p} \wedge \ldots \wedge X_{i_{n} p}
\end{eqnarray*}
because at each point the cup product $\smile$ coincides with tensor product, once each $X_{i p}$ can be viewed as a linear functional. Thus, $\psi(\eta)=0$ results $\eta_{p}=0$ for all $p$, and then $\eta=0$. Lets show that $\psi$ is surjective. Let $N \in Alt(MDer^{n}(\funcman{\manifold{M}}))$. By linearity, it is enough to consider $N$ in the form $\displaystyle \frac{1}{n!}\sum_{\sigma \in S_{n}}\varepsilon(\sigma)X_{\sigma(1)} \smile \ldots \smile X_{\sigma(n)}$. Now, take $\eta \in \Omega_{n}(\manifold{M})$ as $X_{1} \wedge \ldots \wedge X_{n}$. It follows that
\begin{displaymath}
\psi(X_{1} \wedge \ldots \wedge X_{n}) = Alt(X_{1} \smile \ldots \smile X_{n}) = \frac{1}{n!}\sum_{\sigma \in S_{n}}\varepsilon(\sigma)X_{\sigma(1)} \smile \ldots \smile X_{\sigma(n)}
\end{displaymath}
By linearity, $\psi(\eta) = N$. Hence, we have an one-to-one association between alternating elements in $MDer^{n}(\funcman{\manifold{M}})$ and $n$-vector fields. For now on, we shall no more distinguish such elements. We call $J_{n}$ the family of maps taking cochains $D \in D_{poly}^{n,r}(\funcman{\manifold{M}})$ and sending to $J_{n}(D) = Alt(D)$, for $n \geq 1$ and $J_{0}$ as identity on $\funcman{\manifold{M}}$. As $\funcman{\manifold{M}}$ is commutative, $\delta_{H}$ vanish on $\funcman{\manifold{M}}$. Thus, $J_{1} \circ \delta_{H} = d \circ J_{0}$. Let $D$ be a $n$-coboundary, $n > 1$. Then there exists a $(n-1)$-cochain $E$ such that $D=\delta_{H}(E)$. The formulae \ref{formula_dif_hochschild_on_SDer_r} and \ref{formula_dif_hochschild_on_D_poly} shows that $\delta_{H}(E)$ is a linear combination of terms which are symmetric on two entries, hence it must be $Alt(\delta_{H}(E)) = 0$. It follows that $J_{n} \circ \delta_{H} = d \circ J_{n-1}$, because $d$ is identically null. Hence, each $J_{n}$ induces a morphism on cohomology $J_{n}^{\ast}:H^{n}(\mathcal{D}(\funcman{\manifold{M}})) \rightarrow H^{n}(\Omega_{n}(\manifold{M}))$.

By the fact that $d$ is the null differential on $(\Omega_{n}(\manifold{M}),d)$ we have $H^{n}(\Omega_{n}(\manifold{M}))$ isomorphic as $\R$-vector space to $\Omega_{n}(\manifold{M})$, for all $n \geq 0$.

It is clear that $J_{0}^{\ast}$ is isomorphism. Let $D$ be a $n$-cocycle, $n \geq 1$. From proposition \ref{cocycle_coboundary_n_vector_hochschild_formula} we have $D=\delta_{H}(E)+\eta$, where $E$ is a $(n-1)$-cochain and $\eta \in \Omega_{n}(\manifold{M})$. It follows that if $\theta \in H^{n}(D_{poly}(\funcman{\manifold{M}}))$, whose representing element in $D_{poly}(\funcman{\manifold{M}})$ is $D$, then $D$ can be written as $D=\delta_{H}(E)+\eta$ and thus
\begin{displaymath}
J_{n}^{\ast}(\theta)=[J_{n}(D)]=[J_{n}(\delta_{H}(E)+\eta)]=[\eta]=\eta
\end{displaymath}
$J_{n}^{\ast}$ is injective. Indeed, if $\theta$ is such that $J_{n}^{\ast}(\theta) = 0$, then $[J_{n}(D)] = 0$ hence $J_{n}(\delta_{H}(E)+\eta) = J_{n}(\eta) = 0$, resulting $\eta=0$ because $J_{n}(\eta) = \eta$. Thus, $D = \delta_{H}(E)$ and then $\theta$ is the null class. Now, $J_{n}^{\ast}$ is surjective. To show this, note that $\Omega_{n}(\manifold{M})$ is isomorphic to $Alt(MDer^{n}(\funcman{\manifold{M}}))$, which is contained in $MDer^{n}(\funcman{\manifold{M}})$, which is contained in $D_{poly}^{n,r}(\funcman{\manifold{M}})$, for all $r \geq 1$. Hence, given $\eta \in \Omega_{n}(\manifold{M})$ we associate $\eta \in \Omega_{n}(\manifold{M})$ to it. However, by theorem \ref{thm_MDer_trivial_subcomplex}, $\eta$ is a $n$-cocycle. Thus, $J_{n}(\eta)=\eta$. Also, by $\eta$ alternating and by formulae \ref{formula_dif_hochschild_on_SDer_r} and \ref{formula_dif_hochschild_on_D_poly}, $\eta$ can not be a coboundary, therefore the class of $\eta$ in $H^{n}(D_{poly}(\funcman{\manifold{M}}))$ can not be the null class. It follows that $J_{n}^{\ast}$ is an isomorphism on cohomology for all $n$ and hence $(\mathcal{D}(\funcman{\manifold{M}}),\delta_{H})$ and $(\Omega_{\bullet}(\manifold{M}),d)$ are quasi-isomorphics.
\end{proof}\label{hkr}
\end{thm}

\bibliographystyle{plain}
\bibliography{lhp_hkr}

\end{document}